\documentclass[12pt]{article}

\usepackage{latexsym,amsfonts,amsmath,amssymb,graphicx}

\newcommand{\EE}{\mbox{\bf E}\,}
\newcommand{\PP}{\mbox{\bf P}\,}
\newcommand{\CE}{\mbox{CE}}

\newcommand{\de}{\mbox{deg}}
\newcommand{\RW}{\mbox{RW}}
\newcommand{\CRW}{\mbox{CRW}}

\newcommand{\R}{\mathbb{R}}
\newcommand{\C}{\mathbb{C}}

\newcommand{\N}{\mathbb{N}}
\newcommand{\D}{\mathbb{D}}
\newcommand{\Z}{\mathbb{Z}}
\newcommand{\Ree}{\mbox{Re}\,}
\newcommand{\Imm}{\mbox{Im}\,}
\newcommand{\A}{\mbox{\bf A}}
\newcommand{\CC}{\mbox{\bf C}}
\newcommand{\B}{\mbox{\bf B}}
\newcommand{\SA}{\mbox{\bf S}}
\newcommand{\pa}{\partial}
\def\eps{\varepsilon}
\def\til{\widetilde}
\def\ha{\widehat}
\def\sem{\setminus}
\def\lin{\overline}
\def\vphi{\varphi}
\def\up{\upsilon}

\newtheorem{Lemma}{Lemma}[section]
\newtheorem{Theorem}{Theorem}[section]
\newtheorem{Definition}{Definition}[section]
\newtheorem{Corollary}{Corollary}[section]
\newtheorem{Proposition}{Proposition}[section]
\numberwithin{equation}{section}

\begin{document}
\title{Stochastic Loewner evolution \\in doubly connected domains}
\date{March, 2004}
\author{Dapeng Zhan}
\maketitle
\begin{abstract} This paper introduces the annulus SLE$_\kappa$ processes in doubly 
connected domains. 
Annulus SLE$_6$ has the same law as stopped radial SLE$_6$, up to a time-change.
For $\kappa\not=6$, some weak equivalence relation exists between annulus SLE$_\kappa$
and radial SLE$_\kappa$. Annulus SLE$_2$ is the scaling limit of the corresponding
loop-erased conditional random walk, which implies that a certain form of SLE$_2$ 
satisfies the reversibility property. We also consider the disc SLE$_\kappa$ process
defined as a limiting case of the annulus SLE's. Disc SLE$_6$ has the same
law as stopped full plane SLE$_6$, up to a time-change. Disc SLE$_2$ is the 
scaling limit of loop-erased random walk, and is the reversal of radial SLE$_2$. 
\end{abstract}

\section{Introduction}
Stochastic Loewner evolution (SLE), introduced by O.\ Schramm in \cite{14}, is a 
family of random growth processes of plane sets in simply connected domains. 
The evolution is described by the classical
Loewner differential equation with the driving term being a one-dimensional 
Brownian motion. SLE depends on a
parameter $\kappa>0$, the speed of the Brownian motion, and behaves differently
for different value of $\kappa$. See \cite{13} by S.\ Rohde and O.\ Schramm for the basic 
fundamental properties of SLE.

Schramm's processes turned out to be very useful. On the one hand, they are
amenable to computations, on the other hand, they are 
related with some statistical physics models. 
In a series of papers \cite{3}-\cite{7}, G.\ F.\ Lawler, O.\ Schramm and W.\ Werner
used SLE to determine the Brownian motion intersection 
exponents in the plane, identified SLE$_2$ and SLE$_8$ with the scaling limits of
LERW and UST Peano curve, respectively, and conjectured that SLE$_{8/3}$ is the
scaling limit of SAW. S.\ Smirnov proved in \cite{15} that
SLE$_6$ is the scaling limit of critical site percolation on the triangular lattice. 

For various reasons, a similar theory should also exist for multiply connected domains 
and even for general Riemann surfaces. We expect that the definition and some study of 
general SLE will give us better understanding of SLE itself and its physics background. 
The definition of SLE in simply connected domains uses the fact that the complement of 
SLE stopped at a finite time in a simply connected domain other than $\C$ is still simply connected, so it is conformally equivalent to the whole domain. But this property does not 
hold for general domains. That is the main difficulty in our definition of general SLE. 

As a start, we consider SLE in the most simple non-simply connected domains: doubly connected
domains.
We show that the corresponding processes, the annulus 
SLE$_\kappa$, have features similar to those in the simply connected case. More specifically,
we prove that annulus SLE$_6$ has locality property; and for all $\kappa>0$,
annulus SLE$_\kappa$ is equivalent to radial SLE$_\kappa$. We also justify this definition 
by proving that annulus SLE$_2$ is the scaling limit of the corresponding loop-erased 
conditional random walk.

After these, we define disc SLE in simply connected domains, which is the limit
case of annulus SLE. Disc SLE$_6$ also has locality property, so its final hull has the same
law as the hull generated by a plane Brownian motion stopped on hitting the boundary.
Disc SLE$_2$ is the scaling limit of the corresponding loop-erased random walk.
It then follows that disc SLE$_2$ is the
reversal of radial SLE$_2$ started from a random point on the boundary with harmonic measure.  

\subsection{SLE in simply connected domains}
For $\kappa\ge 0$, the standard radial SLE$_\kappa$ is obtained by solving the Loewner differential equations:
$$\partial_t\vphi_t(z)=\vphi_t(z)\frac{1+\vphi_t(z)/\chi_t}{1-\vphi_t(z)/\chi_t}\mbox{, }\mbox{ }
0\le t<\infty\mbox{, }\mbox{ }\mbox{ }\vphi_0(z)=z$$
where $$\chi_t=\exp(iB(\kappa t)),$$ and $B(t)$ is a standard Brownian motion on $\R$
started from $0$. Let $K_t$
be the set of points $z$ in $\D=\{z\in\C:|z|<1\}$ such that the solution $\vphi_s(z)$ blows up before or at time $t$. Then $D_t:=\D\sem K_t$ is a simply connected domain, $0\in D_t$, and $\vphi_t$ maps
$D_t$ conformally onto $\D$ with $\vphi_t(0)=0$ and $\vphi_t'(0)=e^t$. 
The family of hulls $(K_t,0\le t<\infty)$ grows in $\D$ from $1$ to $0$, and is called the
standard radial SLE$_\kappa$. If $\Omega$ is a simply connected domain (other than $\C$), 
$a$ a prime end, $b\in\Omega$, then SLE$_\kappa(\Omega;a\to b)$, radial
SLE$_\kappa$ in $\Omega$ from $a$ to $b$, is defined as the image of the standard 
radial SLE$_\kappa$ under the conformal map $(\D;1,0)\to(\Omega,a,b)$.
By construction, radial SLE is conformally invariant. 

Suppose $(K_t)$ is a radial SLE$_\kappa(\Omega;a\to b)$. 
Then for any fixed $s\ge 0$, the law of a certain conformal image of $(K_{s+t}\sem K_s)$ is 
the same as the law of $(K_t)$, and is independent of $(K_r)_{0\le r\le s}$. 
In other words, radial SLE$_\kappa$ has ``i.i.d." increments, in the sense of conformal 
equivalence. This property, together with the symmetry of the law in $(\D;1,0)$ w.r.t.
complex conjugation, characterizes radial SLE up to $\kappa$. 

Chordal SLE$_\kappa$ processes are defined in a similar way. In this case, the family of 
hulls $(K_t)$ grows in a simply connected domain from one boundary point (prime end) to another.
Once again, the properties of conformal invariance, ``i.i.d." increments, and the
corresponding symmetry property determine a one-parameter family of such processes. 

Radial SLE and chordal SLE are equivalent in the following sense. 
Suppose $\Omega$ is a simply connected
domain, $a$ and $c$ are two distinct prime ends, and $b\in \Omega$. For a fixed $\kappa>0$, let $(K_t)$ be a radial SLE$_\kappa(\Omega;a\to b)$ and $(L_s)$ a chordal SLE$_\kappa(\Omega;a\to c)$. 
Let $T$ be the first time that $K_t$ swallows $c$, $S$ the first time that $L_s$ swallows $b$. 
We set $T$ or $S$ to be $\infty$ by convention if the corresponding hitting time does not exist. 
If $\kappa=6$, up to a time-change, the law of $(K_t)_{0\le t\le T}$ is the same as the law 
of $(L_s)_{0\le s\le S}$. If $\kappa\ne 6$, there exist two sequences of stopping times
$\{T_n\}$ and $\{S_n\}$ such that $T=\lor_n T_n$, $S=\lor_n S_n$, and for each $n\in\N$,
up to a time-change, the laws of $(K_t)_{0\le t\le T_n}$ and $(L_s)_{0\le s\le S_n}$ are
equivalent. In other words, they have positive density w.r.t.\ each other. 
The strong equivalence relation of radial and chordal
SLE$_6$ is related to the so-called locality property: the SLE$_6$ hulls do not feel the
boundary before hitting it.

The equivalence property ensures that 
for the same $\kappa$, radial SLE$_\kappa$ and chordal SLE$_\kappa$ behave similarly. 
For instance, if $\kappa\le 4$, and $(K_t)$ is a radial or chordal 
SLE$_\kappa$ in $\Omega$, then a.s.\ there is a simple path 
$\beta:(0,\infty)\to \Omega$ such that for any $t\in[0,\infty)$, we have $K_t=\beta(0,t]$. 
If $\kappa>4$ and $\partial \Omega$ is locally connected, then 
a.s.\ there is a non-simple path $\beta:(0,\infty)\to\lin{\Omega}$ such that for any $t\in[0,\infty)$, $K_t$ is the hull generated by $\beta(0,t]$. 
This path $\beta$ is called the SLE$_\kappa$ trace.  

Full plane SLE$_\kappa$: $(K_t,-\infty<t<\infty)$ grows in $\ha{\C}=\C\cup\{\infty\}$ 
from $0$ to $\infty$. For any fixed $s\in\R$, the law of a certain conformal image of
$(K_{s+t}\sem K_s)$ is the same as the law of 
the standard radial SLE$_\kappa$, and is independent of $(K_r)_{-\infty
<r\le s}$. Full plane SLE can be viewed as the limit of radial 
SLE$_\kappa(\ha{\C}\sem\eps\lin{\D};\eps\to\infty)$ as $\eps\to 0^+$. 

\subsection{Definition of annulus SLE}

For $p>0$, we denote by $\A_p$ the standard annulus of modulus $p$:
$$\A_p=\{z\in\C:e^{-p}<|z|<1\}.$$ Every doubly connected domain $D$ with non-degenerate 
boundary is conformally equivalent to a unique $\A_p$, and $p=M(D)$ is the modulus 
of $D$. We may first define SLE on the standard annuli, and then extend the definition to 
arbitrary doubly connected domains via conformal maps. 

Denote $$\SA_p(z)=\lim_{N\to\infty}\sum_{-N}^N\frac{e^{2kp}+z}{e^{2kp}-z}.$$
For $\chi\in\partial D$, let 
$$\SA_p(\chi,z)=\SA_p(z/\chi).$$ 
The function $\SA_p(\chi,\cdot)$ is a Schwarz kernel of $\A_p$ in the sense that 
if $f$ is an analytic function in $\A_p$, continuous up to the boundary, 
and constant on the circle $\CC_p:=\{z\in\C:|z|=e^{-p}\}$, then for any $z\in\A_p$,
 $$f(z)=\int_{\CC_0}f(\chi)\SA_p(\chi,z)d{\bf m}+iC,$$
where $\bf m$ is the uniform probability measure on $\CC_0=\partial\D$, and
$C$ is some real constant. Note that the Schwarz kernels are not unique. The choice of 
$\SA_p(\chi,\cdot)$ here satisfies the rotation symmetry and reflection symmetry.

Let $\chi:[0,p)\to\CC_0$ be a continuous function. Consider the following 
Loewner-type differential equation: \begin{equation}
\partial_t\vphi_t(z)=\vphi_t(z)\SA_{p-t}(\chi_t,\vphi_t(z))\mbox{, }\mbox{ }
0\le t<p\mbox{, }\mbox{ }\mbox{ }\vphi_0(z)=z.
\end{equation}
For $0\le t<p$, let $K_t$ be the set of $z\in\A_p$ such that the solution $\vphi_s(z)$ blows 
up before or at time $t$. Let $D_t=\A_p\sem K_t$, $0\le t<p$. We call $K_t$ ($\vphi_t$, resp.), $0\le t<p$, the standard annulus LE hulls (maps, resp.) of modulus $p$ driven 
by $\chi_t$, $0\le t<p$. We will see that for each $0\le t<p$, $\vphi_t$ maps $D_t$ conformally
onto $\A_{p-t}$, and maps $\CC_p$ onto $\CC_{p-t}$. 

If we replace $\SA_{p-t}(\chi_t,\vphi_t(z))$ in formula (1.1) by
$$\ha\SA_{p-t}(\chi_t,\vphi_t(z)):=\SA_{p-t}(\chi_t,\vphi_t(z))-
\Imm \SA_{p-t}(\chi_t,e^{t-p}),$$
and let $\ha\vphi_t(z)$ be the corresponding solutions. Then we have $\ha\vphi_t(e^{-p})
=e^{t-p}$, $0\le t<p$, since $\ha\SA_{p-t}(\chi_t,e^{t-p})\equiv 1$. 
Actually $\ha\SA_p$ is the Schwarz kernel in \cite{16}. We will use it in the
proof of Proposition 2.1. We prefer $\SA_p$ to $\ha\SA_p$ in the definition of annulus 
SLE because 
if we use $\ha\SA_p$ then the driving function must contain a drift term besides a Brownian
motion. See the definition of SLE$_6$ in \cite{b}.

We define standard annulus SLE$_\kappa$ of modulus $p$ to be the solution of (1.1) with 
$\chi_t=\exp(iB(\kappa t))$, $0\le t<p$. The family of hulls grows from $1$ to $\CC_p$. 
Via a certain
conformal map, we may extend the definition to SLE$_\kappa(\Omega;a\to B)$ where $\Omega$
is a doubly connected domain with non-degenerate boundary, $B$ is a boundary component,
and $a$ is a boundary point (prime end) on the other boundary component.
Note that the conformal type of $\Omega\sem K_t$ is always changing, so the annulus
SLE$_\kappa$ hulls cannot have identical increments in the sense of conformal equivalence. 
We may only require that
for any fixed $s\in[0,p)$, the conformal image of $(K_{s+t}\sem K_s)_{0\le t<p-s}$
has the same law as the annulus SLE hulls of modulus $p-s$. This together with 
the symmetry property does not determine the driving process up to a single parameter. 
However, it turns out that $\exp(iB(\kappa t))$, a Brownian motion on $\CC_0$ started from 
$1$ with constant speed $\kappa$, is a reasonable choice for the driving 
process. The main goal of the paper is to justify this claim.

Two facts of doubly connected domains are used in the above definition of annulus SLE.
First, the conformal type of a doubly connected domain can be described by a single number,
which is the modulus. So we use the time parameter to describe the modulus. 
Second, given a boundary component $B$ and a prime end $P$ on the other boundary component of
some doubly connected domain $D$, there is a self-conjugate-conformal map of $(D;B,P)$. 
This is clear when $D$ is the standard annulus.
We actually assume that the law of annulus SLE$_\kappa(D;P\to B)$ is invariant under that map. Because of these, our definition of annulus SLE can
be expressed by some nice differential equations. However, these two facts do not hold
for $n$-connected domains when $n>2$. Some other methods are needed to define the SLEs. 
The extensions of SLE to multiply connected domains and Riemann surfaces are now in preparation,
and will appear elsewhere.

\subsection{Main results}
Suppose $\Omega$ is a simply connected domain, $a$ is a prime end, and $b$ is an interior point. Suppose $F
\supsetneqq\{b\}$ is a contractible compact subset of $\Omega$. Then $\Omega\sem F$ is 
a doubly connected domain with two boundary components $\partial \Omega$ and $\partial F$. 
We call $F$ a hull in $\Omega$ w.r.t.\ $b$. 
For a fixed $\kappa>0$, let $(K_t)$ be a radial SLE$_\kappa(\Omega;a\to b)$, and $(L_s)$ an 
annulus SLE$_\kappa(\Omega\sem F;a\to\partial F)$. Then we have

\begin{Theorem} (i)If $\kappa=6$, the law of $(K_t)_{0\le t<T_F}$, is equal to that of $(L_s)_{0\le s<p}$, up to a time-change. \\
(ii)If $\kappa\not=6$, there exist two sequences of stopping times $\{T_n\}$ and 
$\{S_n\}$ such that $T=\lor_n T_n$, $p=\lor_n S_n$, and for each $n\in\N$,
the law of $(K_t)_{0\le t\le T_n}$ is equivalent to that of $(L_s)_{0\le s\le S_n}$, 
up to a time-change.
\end{Theorem}

The second main result of the paper concerns the convergence of a loop-erased conditional 
random walk (LERW) with appropriate boundary conditions to an annulus SLE$_2$.
For any plane domain $\Omega$, and $\delta>0$, let $\Omega^\delta$ denote the graph defined as follows. The vertex set $V(\Omega^\delta)$ consists of the points in $\delta\Z^2\cap \Omega$ 
and the intersection points of $\partial \Omega$ with edges of $\delta\Z^2$. The edge set
$E(\Omega^\delta)$ consists of 
the unordered vertex pairs $\{u,v\}$ such that the line segment $(u,v)\subset \Omega$, and 
there is an edge of $\delta\Z^2$ that contains $(u,v)$ as a subset.

Suppose $D$ is a doubly connected domain with boundary components $B_1$ and $B_2$, 
$0\in B_1$ and there is some $a>0$ such that the line segment $(0,a]$ is contained in $D$. 
This line segment determines a prime end in $D$ on $B_1$, denoted by $0_+$.
We may assume that $\delta$ is sufficiently small so that
$0$ and $\delta$ are adjacent vertices of $D^\delta$, and there is a lattice path 
on $D^\delta$ connecting $\delta$ and $V(D^\delta)\cap B_2$. 

Now let $\RW$ be a simple random walk on $D^\delta$ started from $\delta$ and stopped on hitting 
$\partial D$. Let $\CRW$ be $\RW$ conditioned on the event that $\RW$ hits $B_2$ before $B_1$. 
Let LERW be the loop-erasure of CRW, which
is obtained by erasing the loops of $\CRW$ in the order that they appear. See \cite{2} for details. 
Then LERW is a random simple lattice path on $D^\delta$ from $\delta$ to $B_2$.
We may also view LERW as a random simple curve in $D$ from $\delta$ to $B_2$. 
Taking with the segment $[0,\delta]$, we obtain a random simple curve in $D$ from $0$ to $B_2$. We parameterize this curve by $\beta^\delta[0,p]$ so that $\beta^\delta(0)=0$, 
$\beta^\delta(p)\in B_2$, and $M(D\sem\beta^\delta(0,t))=p-t$, for $0\le t<p$. 

Now let $(K^0_t)_{0\le t<p}$ be an annulus SLE$_2(D;0_+\to B_2)$. 
From Theorem 1.1 and the existence of radial SLE$_\kappa$ traces, we know that a.s.\ there 
exists a random simple path $\beta^0(t)$, $0<t<p$, such that 
$K^0_t=\beta^0(0,t]$, for $0\le t<p$. 

\begin{Theorem} For every $q\in(0,p)$ and $\eps>0$, there is a $\delta_0>0$ depending on $q$ and $\eps$ such that 
for $\delta\in(0,\delta_0)$ there is a coupling of the processes $\beta^\delta$ and $\beta^0$ 
such that 
$$\mbox{\bf P}[\sup\{|\beta^\delta(t)-\beta^0(t)|:t\in[q,p)\}>\eps]<\eps.$$
Moreover, if the impression of the prime end $0_+$ is a single point, then the theorem holds 
with $q=0$.
\end{Theorem}

Here a coupling of two random processes $A$ and $B$ is a probability space with two random processes $A'$ and $B'$,
where $A'$ and $B'$ have the same law as $A$ and $B$, respectively. In the above statement 
(as is customary)
we don't distinguish between $A$ and $A'$ and between $B$ and $B'$. The impression (see \cite{11}) of a prime
end is the intersection of the closure of all neighborhoods of that prime end. 

For $\kappa=2$, $8$ and $8/3$, chordal SLE$_\kappa$ satisfies the 
reversibility property. That means the reversal of chordal SLE$_\kappa(D;a\to b)$ trace 
has the same law as chordal SLE$_\kappa(D;b\to a)$ trace, up to a time-change. 
For the annulus SLE trace, the starting point is a fixed prime end, but the end point 
(if it exists)
is a random point on a boundary component. To get the reversibility property, we 
have to ``average" the annulus SLE traces in the same domain started from different points 
of one boundary component. From Theorem 1.2 and the reversibility of LERW (see \cite{2}), it 
then follows

\begin{Corollary} The reversal of the annulus SLE$_2(\A_p;{\bf x}\to\CC_p)$ trace has the 
same law as the annulus SLE$_2(\A_p;{\bf y}\to\CC_0)$ trace, up to a time-change,
where $\bf x$ and $\bf y$ are uniform random points on $\CC_0$ and $\CC_p$, respectively.\end{Corollary}

The definition of annulus SLE enables us to define disc SLE$_\kappa$ that grows in a simply 
connected domain $\Omega$ from an interior point to 
the whole boundary. It can be viewed as the limit of annulus SLE$_\kappa$ as the modulus 
tends to infinity. The relation between disc SLE and annulus SLE is similar to that between 
full plane SLE and radial SLE.

From our methods, it follows that for any simply connected domain $\Omega$ that contains $0$, 
the full plane SLE$_6$ before the hitting time of $\partial \Omega$ has the same law as 
the disc SLE$_6(\Omega;0\to\partial \Omega)$, up to a time-change. 
This gives an alternative proof of the following facts mentioned in \cite{17}\cite{7}. 
The hitting point of full plane SLE$_6$ at $\partial \Omega$ has harmonic measure valued at
$0$, and therefore the full plane SLE$_6$ hull at the hitting time of $\partial \Omega$ has
the same law as the hull generated by a plane Brownian motion started from $0$ and stopped
on exiting $\Omega$.

We also show that the LERW on the grid approximation $\Omega^\delta$ started from an interior 
vertex $0$ to the boundary converges to the disc SLE$_2(\Omega;0\to\partial \Omega)$, 
as $\delta\to 0$. Together with the approximation result in \cite{6}, this implies that 
the reversal of the disc SLE$_2(\Omega;0\to\partial \Omega)$ has the same law as the 
radial SLE$_2(\Omega;{\bf z}\to 0)$, up to a time-change, where $\bf z$ is a random point
on $\partial \Omega$ that has harmonic measure valued at $0$. 

\subsection{Some comments about the proof}
The discussion of the convergence of LERW to annulus SLE$_2$ basically follows the methods
developed in \cite{6}. In the same order as in \cite{6}, logically, we first find the 
observables for LERW; then prove they are martingales and converge to some continuous 
harmonic functions; these facts are used to show that the driving function of the LERW 
converges to the Brownian motion with speed $2$; finally we use the nice behavior of LERW 
path to show that the path parameterized according to the modulus of the remaining domain
converges to the annulus SLE$_2$ trace uniformly in probability. 

However, some notations and proofs in \cite{6} can not be transplanted
to this paper immediately. For example, the observables in this paper has counterparts in
simply connected domains, which are exactly the observables introduced in \cite{6}. 
But the LERW studied there is from an interior vertex to the boundary, and the proof 
of Proposition 3.4 in \cite{6} uses this construction. We have to prove the fact
that they are martingales using a different method, which we believe shows 
some essence of this subject. Moreover, since the
moduli change in time, some proofs here, e.g., that to Proposition 3.4, are much longer than their counterparts, e.g., part of the proof to Proposition 3.4 in \cite{6}.

The authors of \cite{6} first use some subgraph of $\Z^2$ to approximate a simply connected
plane domain, and they use the inner radius with respect to a fixed point (which is $0$ there) 
to describe the extent that the graph approximates the domain. 
After some rescaling, the inner radius means the distance from $0$ to the boundary of the 
domain divided by the length of the mesh. It seems not easy to 
find counterparts of the inner radius for doubly connected domains. So we proceed in another 
way by taking the limit of some sequence of domains. This results in a very long proof of Proposition 3.3. This method extends to the cases of multiply connected domains.

\section{Equivalence of annulus and radial SLE}
\subsection{Deterministic annulus LE hulls}
We recall some facts about the Schwarz function 
$$\SA_r(z)=\lim_{N\to\infty}\sum_{k=-N}^N\frac{e^{2kr}+z}{e^{2kr}-z}\mbox{, }\mbox{ }r>0.$$ 
(i) $\SA_r$ is analytic in $\C\sem\{0\}\sem\{e^{2kr}:k\in\Z\}$;\\
(ii) $\{e^{2kr}:k\in\Z\}$ are simple poles of $S_r$;\\ 
(iii) $\Ree \SA_r\equiv 1$ on $\CC_r=\{z\in\C:|z|=e^{-r}\}$;\\
(iv) $\Ree \SA_r\equiv 0$ on $\CC_0\sem\{1\}$;\\
(v) $\Ree \SA_r>0$ in $\A_r$; and \\
(vi) $\Imm \SA_r\equiv 0$ on $\R\sem\{0\}\sem\{\mbox{poles}\}$. 

Suppose $f$ is an analytic function in $\A_r$, $\Ree f$ is non-negative, and $\Ree f(z)$ tends to $a$ as 
$z\to \CC_r$, then there is some positive measure $\mu=\mu(f)$ on $\CC_0$ of total mass $a$
 such that 
\begin{equation} f(z)=\int_{\CC_0}\SA_r(z/\chi)d\mu(\chi)+iC,\end{equation}
for some real constant $C$. If $\Ree f(z)$ tends to zero as $z$ approaches the complement of an 
arc $\alpha$ of $\CC_0$, then $\mu(f)$ is supported by $\lin{\alpha}$. 
Moreover, if $f$ is bounded, then the radial limit of $f$ on $\CC_0$ exists a.e.,
and $d\mu(f)/d{\bf m}=f|_{C_0}$. The proof is similar to that of the Poisson integral formula. 

Divide both sides of equation (1.1) by $\vphi_t(z)$and take the real part. We get
$$\partial_t\ln|\vphi_t(z)|=\Ree \SA_{p-t}(\vphi_t(z)/\chi_t).$$
From the values of $\Ree \SA_{p-t}$ on $\CC_{p-t}$ and $\CC_0$ we see that 
if $z\in\CC_0\sem\{1\}$, then 
$\vphi_t(z)\in\CC_0\sem\{1\}$ until it blows up; if $z\in\CC_p$, then $\vphi_t(z)\in\CC_{p-t}$ for $0\le t<p$. 
Thus for $z\in\A_p$, $\vphi_t(z)$ stays between $\CC_0$ and $\CC_{p-t}$ until it blows up. So $\vphi_t$ maps
$D_t$ into $\A_{p-t}$. The fact that $\SA_{p-t}$ is analytic implies that for every $t\in[0,p)$, $\vphi_t$ is a conformal map of $D_t$. 
By considering the backward flow, it is easy to see that 
$\vphi_t$ maps $D_t$ onto $\A_{p-t}$.

\begin{Definition} 
Suppose $D$ is a doubly connected domain with boundary components $B$ and $B'$. 
We call $K\subset D$ a hull in $D$ on $B$ if $D\sem K$ is a doubly connected domain that has $B'$ as a boundary component. The capacity of $K$ in $D$ w.r.t.\ $B'$, denoted by $C_{D,B'}(K)$,
 is the value of $M(D)-M(D\sem K)$. \end{Definition}

\begin{Definition} Suppose $\Omega$ is a simply connected domain. We call
$K\subset \Omega$ a hull in $\Omega$ on $\partial \Omega$, if $\Omega\sem K$ is a simply connected domain.
If $\vphi$ maps $\Omega\sem K$ conformally onto $\Omega$ and for some $a\in \Omega\sem K$, 
$\vphi(a)=a$ and $\vphi'(a)>0$, then $\ln\vphi'(a)>0$, and is called the capacity of 
$K$ in $\Omega$ w.r.t.\ $a$, denoted by $C_{\Omega,a}(K)$. \end{Definition}

If $K$ is a hull in $\A_p$ on $\CC_0$, and 
$\psi$ is any conformal map from $\A_p\sem K$ onto $\A_{p-r}$ which takes $\CC_p$ to 
$\CC_{p-r}$, then the radial limit of $\psi^{-1}$ on $\CC_0$ exists a.e.,
and $$C_{A_p,C_p}(K)=\int_{\CC_0}-\ln|\psi^{-1}|d{\bf m}.$$

If $K$ is a hull in $\D$ on $\CC_0$ and $\vphi$ maps $\D\sem K$ onto $\D$ conformally 
so that $\vphi(0)=0$, then the radial limit of $\vphi^{-1}$ on $\CC_0$ exists 
a.e., and $$C_{\D,0}(K)=\int_{\CC_0}-\ln|\vphi^{-1}|d{\bf m}.$$




Similarly as Lemma 2.8 in \cite{3}, using the 
integral formulas for capacities of hulls in $\D$ and $\A_p$, it is not hard to derive 
the following Lemma:
\begin{Lemma} Suppose $x,y\in\CC_0$, and $G$ is a conformal map from a neighborhood
$U$ of $x$ onto a neighborhood $V$ of $y$ such that $G(U\cap\D)=V\cap\D$.
Fix any $p>0$. For every $\eps>0$, there is $r=r(\eps)>0$ such that if  $K$ is 
a non-empty hull in $\D$ on $\CC_0$ and $K\subset\B(x;r)$, the open ball of radius $r$ about 
$x$, then $K\subset U$, $G(K)$ is 
a hull in $\A_p$ on $\CC_0$, and 
$$\left|\frac{C_{A_p,C_p}(G(K))}{C_{\D,0}(K)}-|G'(x)|^2\right|<\eps.$$
\end{Lemma}

Suppose $D$ is a doubly connected domain with boundary components $B_1$ and $B_2$. 
We call $(K_s, a\le s<b)$ a Loewner chain in $D$ on $B_1$ if every $K_s$ is a 
hull in $D$ on $B_1$, $K_{s_1}\subsetneqq K_{s_2}$ if $a\le s_1<s_2<b$, and for every $c\in(a,b)$,
the extremal length (see \cite{1}) of the family of curves in $D\sem K_{s+u}$ that disconnect
 $K_{s+u}\sem K_s$ from $B_2$ tends to $0$ as $u\to 0^+$, uniformly in $s\in[a,c]$.
If the area of $D$ is finite, then the above condition holds iff 
the infimum length of all $C^1$ curves in $D\sem K_{s+u}$ that disconnect
$B_2$ from $K_{s+u}\sem K_s$ tends to $0$ as $u\to 0^+$, uniformly in $s\in[a,c]$.

Now we consider a Loewner chain in $\A_p$ on $\CC_0$. The following proposition is
similar to the theorems for chordal and radial LE in \cite{3} and \cite{9}.

\begin{Proposition} The following two statements are equivalent:\\
1. $K_t$, $0\le t<p$, are the standard LE hulls of modulus $p$ driven by some continuous function
 $\chi:[0,p)\to\CC_0$;\\
2. $(K_t,0\le t<p)$ is a Loewner chain in $\A_p$ on $\CC_0$, and  $C_{A_p,C_p}(K_t)=M(\A_p)-M(\A_p\sem K_t)=t$ for $0\le t<p$.\\
Moreover, $\{\chi_t\}=\cap_{u>0}\lin{\vphi_t(K_{t+u}\sem K_t)}$, where $\vphi_t$ is the 
standard annulus LE map. If $(L_s,a\le s<b)$ is any Loewner chain 
in $\A_p$ on $\CC_0$, then $s\mapsto C_{A_p,C_p}(L_s)$ is a continuous 
(strictly) increasing function.
\end{Proposition}
{\bf Proof.} The method of the proof is a combination of extremal length comparison,
the use of formula (2.1), and some estimation of Schwarz kernels.
It is very similar to the proof of the counterparts in \cite{3} and \cite{9}. 
So we omit the most part of it. One thing we want to show here is how we derive 
$\vphi_t$ from $K_t$ in the proof of 2 implies 1.
We first choose $\ha{\vphi}_t$ that maps $\A_p\sem K_t$ conformally onto $\A_{p-t}$
such that $\ha{\vphi}_t(\CC_p)=\CC_{p-t}$ and $\ha{\vphi}_t(e^{-p})=e^{t-p}$. 
Then we prove that $\ha{\vphi}_t$ satisfies the equation
$$\partial_t\ha{\vphi}_t(z)=\ha{\vphi}_t(z)(\SA_{p-t}(\ha{\vphi}_t(z)/\ha{\chi}_t)-i\Imm 
\SA_{p-t}(e^{t-p}/\ha{\chi}_t)),$$
for some continuous $\ha{\chi}:[0,p)\to\CC_0$. And $\{\ha{\chi}_t\}
=\cap_{u>0}\lin{\ha{\vphi}_t(K_{t+u}\sem K_t)}$.
Define $$\theta(t)=\int_0^t \Imm \SA_{p-s}(e^{s-p}/\ha{\chi}_s)ds,$$ 
$\chi_t=e^{i\theta(t)}\ha{\chi}_t$
and $\vphi_t(z)=e^{i\theta(t)}\ha{\vphi}_t(z)$, for $t\in[0,p)$. Then $\vphi_0(z)=\ha{\vphi}_0(z)=z$, $\vphi_t$ 
maps $\A_p\sem K_t$ conformally onto $\A_{p-t}$, $\{{\chi}_t\}
=\cap_{u>0}\lin{{\vphi}_t(K_{t+u}\sem K_t)}$, and 
$$\partial_t\ln \vphi_t(z)=\partial_t\ln\ha{\vphi}_t(z)+i\theta'(t)=\SA_{p-t}(\ha{\vphi}_t(z)/\ha{\chi}_t)
=\SA_{p-t}(\vphi_t(z)/\chi_t).$$
Thus $\partial_t \vphi_t(z)=\vphi_t(z) \SA_{p-t}(\vphi(z)/\chi_t)$. So $K_t$, $0\le t<p$, are the 
standard annulus LE hulls of modulus $p$, driven by $\chi_t$, $0\le t<p$. $\Box$ 

\subsection{Proof of Theorem 1.1}
We may assume in Theorem 1.1 that $\Omega=\D$, $a=1$ and $b=0$. 
Then $(K_t,0\le t<\infty)$ is the standard radial SLE$_\kappa$. 
Suppose $\vphi_t$ and $\chi_t$, $0\le t<\infty$, are the corresponding standard radial SLE$_\kappa$ maps and 
driving process, respectively. Then $\chi_t=e^{iB(\kappa t)}$, where $B(t)$ is a standard  Brownian motion on $\R$ started from $0$.

For $0\le t<T_F$, $\D\sem F\sem K_t$ is a doubly connected domain. So $K_t$, $0\le t<T_F$, 
are hulls in $\D\sem F$ on $\CC_0$. From \cite{9} we know that $(K_t,0\le t<T_F)$ is a Loewner chain in $\D\sem F$ on $\CC_0$. Suppose $W$ maps $\D\sem F$ conformally onto $\A_p$ so that $W(1)=1$. Then $(W(K_t),0\le t<T_F)$ is a Loewner chain in $\A_p$ on 
$\CC_0$. From \cite{13} we know that $K_t$ approaches $F$ as $t\nearrow T_F$, so $W(K_t)$ 
approaches $\CC_p$ as $t\nearrow T_F$. This implies that $M(\D\sem F\sem K_t)\to 0$ as $t\nearrow T_F$.
Let $u(t)=C_{D,\partial F}(K)=C_{A_p,C_p}(W(K))$. Then 
$u$ is a continuous increasing function and maps $[0,T_F)$ onto $[0,p)$. 
Let $v$ be the inverse of $u$. By Proposition 2.1, $W(K_{v(s)})$, 
$0\le s<p$, are the standard annulus LE hulls of modulus $p$ driven by some continuous 
$\nu: [0,p)\to\CC_0$. 
Let $\psi_s$, $0\le s<p$, be the corresponding standard annulus LE maps. 

Now $\vphi_t$ maps $\D\sem F\sem K_t$ conformally onto $\D\sem\vphi_t(F)$.
Let $f_t=\psi_{u(t)}\circ W\circ\vphi_t^{-1}$. Then $f_t$
maps $\D\sem\vphi_t(F)$ conformally onto $\A_{p-u(t)}$,
and $f_t(\CC_0)=\CC_0$. By Schwarz reflection, we may extend $f_t$ analytically to $\Sigma_t$,
which is the union of $\D\sem\vphi_t(F)$, $\CC_0$, and the reflection of 
$\D\sem\vphi_t(F)$ w.r.t.\ $\CC_0$. And $f_t$ is a conformal map on $\Sigma_t$.  
Note that $f_t$ maps $\vphi_t(K_{t+a}\sem K_t)$ to $\psi_{u(t)}
(W(K_{t+a})\sem W(K_{t}))$ for $a>0$. From Proposition 2.1, we see that 
$\{\nu_{u(t)}\}=\cap_{a>0}\lin{\psi_{u(t)}(W(K_{t+a})\sem W(K_{t}))}$. And from 
the counterpart in \cite{9} of Proposition 2.1, we know that 
$\{\chi_t\}=\cap_{a>0}\lin{\vphi_t(K_{t+a}\sem K_t)}$. Thus $\nu_{u(t)}=f_t(\chi_t)$.
Now $\vphi_t(K_{t+a}\sem K_t)$ is a hull in $\D$, $\vphi_{t+a}\circ\vphi_t^{-1}$ 
maps $\D\sem\vphi_t(K_{t+a}\sem K_t)$ conformally onto $\D$, fixes $0$, and 
$(\vphi_{t+a}\circ\vphi_t^{-1})'(0)=e^a$. So the capacity w.r.t.\ $0$ of $\vphi_t(K_{t+a}\sem K_t)$ is $a$. 
Similarly, $\psi_{u(t)}(W(K_{t+a}\sem W(K_{t}))$ is
a hull in $\A_{p-u(t)}$ on $\CC_0$, and the capacity is $u(t+a)-u(t)$. From Lemma (2.1) we 
conclude that $u'_+(t)=|f_t'(\chi_t)|^2$. 

Let $H=\{(t,z):0\le t<T_F,z\in\Sigma_t\}$ and $G(\chi)=\{(t,\chi_t):0\le t<T_F\}$. 
By the definition of $f_t$, we see that $(t,z)\mapsto f_t'(z)$ is continuous in $H\sem G(\chi)$. 
Note that $f_t'$ is analytic in $\Sigma_t$ for each $t\in[0,T_F)$. The maximum principle
implies that $(t,z)\mapsto f_t'(z)$ is continuous in $H$.
In particular, $t\mapsto f_t'(\chi_t)$ is continuous. So we have

\begin{Lemma} $u(t)$ is $C^1$ continuous, and $u'(t)=|f_t'(\chi_t)|^2$.\end{Lemma}

The fact $W(\chi_0)=W(1)=1$ implies that $\nu_0=1$. We now lift $f_t$ to the covering space. 
Write $\chi_t=e^{i\xi_t}$ and $\nu_s=e^{i\eta_s}$, where $\xi_t=B(\kappa t)$, 
$0\le t<\infty$, and $\eta_s$, $0\le s<p$, is a real continuous function
with $\eta_0=0$. Let $\til{\Sigma}_t=\{z\in\C:e^{iz}\in\Sigma_t\}$. Then there is a unique conformal map $\til{f}_t$ 
on $\til{\Sigma}_t$ such that $e^{i\til{f}_t(z)}=f_t(e^{iz})$ and $\eta_{u(t)}=\til{f}_t(\xi_t)$. And $\til{f}_t$ takes real values on the real line. Moreover, $u'(t)=|f_t'(\chi_t)|^2=\til{f}_t'(\xi_t)^2$. 

\begin{Lemma} $(t,x)\mapsto \til{f}_t(x)$ is $C^{1,\infty}$ continuous on $[0,T_F)\times\R$. And for all $t\in[0,T_F)$, $\partial_t\til{f}_t(\xi_t)=-3\til{f}_t''(\xi_t)$. \end{Lemma}
{\bf Proof.} For any $t\in[0,T_F)$, and $z\in\D\sem F\sem K_t$, we have $f_t\circ\vphi_t(z)=\psi_{u(t)}\circ W(z)$. 
Taking the derivative w.r.t.\ $t$, we compute 
$$\partial_t f_t(\vphi_t(z))+f_t'(\vphi_t(z))\vphi_t(z)\frac{\chi_t+\vphi_t(z)}{\chi_t-\vphi_t(z)}
=u'(t)\psi_{u(t)}(W(z))\SA_{p-u(t)}(\psi_{u(t)}(W(z))/\eta_{u(t)}).$$
By Lemma 2.2, $u'(t)=|f_t'(\chi_t)|^2$. Thus for any $t\in[0,T_F)$ and $z\in\D\sem F\sem K_t$,
$$\partial_t f_t(\vphi_t(z))=|f_t'(\chi_t)|^2f_t(\vphi_t(z))\SA_{p-u(t)}(f_t(\vphi_t(z))
/f_t(\chi_t))-f_t'(\vphi_t(z))\vphi_t(z)\frac{\chi_t+\vphi_t(z)}{\chi_t-\vphi_t(z)}.$$
For any $t\in[0,T_F)$, and $w\in\D\sem\vphi_t(F)$, we have $\vphi_t^{-1}(w)\in\D\sem F\sem K_t$. Thus
$$\partial_t f_t(w)=|f_t'(\chi_t)|^2f_t(w)\SA_{p-u(t)}(f_t(w)/f_t(\chi_t))-f_t'(w)w\frac{\chi_t+w}{\chi_t-w}.$$
Let $g_t(w)$ be the right-hand side of the above formula for $t\in[0,T_F)$ and $w\in\Sigma_t\sem\{\chi_t\}$. 
Then for each $t\in[0,T_F)$, $g_t(w)$ is analytic in $\Sigma_t\sem\{\chi_t\}$.
And $(t,w)\mapsto g_t(w)$ is $C^{0,\infty}$ continuous on $H\sem G(\chi)$.

Now fix $t_0\in[0,T_F)$. Let us compute the limit of $g_{t_0}(w)$ when $w\to\chi_{t_0}$. Since 
$$\SA_{p-u(t_0)}(f_{t_0}(w)/f_{t_0}(\chi_{t_0}))-\frac{f_{t_0}(\chi_{t_0})+f_{t_0}(w)}
{f_{t_0}(\chi_{t_0})-f_{t_0}(w)}\to 0
\mbox{, }\mbox{ as }w\to\chi_{t_0},$$ 
so the limit of $g_{t_0}(w)$ is equal to the limit of the following function:
$$|f_{t_0}'(\chi_{t_0})|^2f_{t_0}(w)\frac{f_{t_0}(\chi_{t_0})+f_{t_0}(w)}{f_{t_0}(\chi_{t_0})-f_{t_0}(w)}-f_{t_0}'(w)w\frac{\chi_{t_0}+w}{\chi_{t_0}-w}.$$  
Let $w=e^{ix}$, we may express the above formula in term of $x$, $\xi_{t_0}$ and $\til{f}_{t_0}$, which is  
$$\til{f}_{t_0}'(\xi_{t_0})^2e^{i\til{f}_{t_0}(x)}\frac{e^{i\til{f}_{t_0}(\xi_{t_0})}+e^{i\til{f}_{t_0}(x)}}
{e^{i\til{f}_{t_0}(\xi_{t_0})}-e^{i\til{f}_{t_0}(x)}}-\til{f}_{t_0}'(x)e^{i\til{f}_{t_0}(x)}
\frac{e^{i\xi_{t_0}}+e^{ix}}{e^{i\xi_{t_0}}-e^{ix}}$$
$$=-ie^{i\til{f}_{t_0}(x)}[\til{f}_{t_0}'(\xi_{t_0})^2\cot(\frac{\til{f}_{t_0}(x)-\til{f}_{t_0}(\xi_{t_0})}{2})
-\til{f}_{t_0}'(x)\cot(\frac{x-\xi_{t_0}}{2})].$$

By expanding the Laurent series of $\cot(z)$ near $0$, we see that the limit of the above formula is $3ie^{i\til{f}_{t_0}(\xi_{t_0})}\til{f}_{t_0}''(\xi_{t_0})=3if_{t_0}(\chi_{t_0})\til{f}_{t_0}''(\xi_{t_0})$. 
Therefore $g_t$ has an analytic extension to $\Sigma_t$ for each $t\in[0,T_F)$. 
The maximum principle also implies that $g_t(w)$ is $C^{0,\infty}$ continuous in $H$, 
and $\partial_t f_t(w)=g_t(w)$ holds in the whole $H$. 
Thus $f_t(w)$ is $C^{1,\infty}$ continuous on
$[0,T_F)\times\CC_0$, and $\til{f}_t(w)$ is $C^{1,\infty}$ continuous on $[0,T_F)\times\R$. 
Finally, $$\partial_t\til{f}_t(\xi_t)=\frac{i\partial_t f_t(\chi_t)}{f_t(\chi_t)}
=\frac{ig_t(\chi_t)}{f_t(\chi_t)}=\frac{-3f_t(\chi_t)\til{f}_t''(\xi_t)}{f_t(\chi_t)}
=-3\til{f}_t''(\xi_t).\mbox{ }\mbox{ }\Box$$
\vskip 2mm
\noindent{\bf Proof of Theorem 1.1.} Note that $\eta_{u(t)}=\til{f}_t(\xi_t)$, 
$\xi_t=B(\kappa t)$, and from Lemma 2.2,
$\partial_t\til{f}_t(\xi_t)=-3\til{f}_t''(\xi_t)$. By It\^o's formula, we have
$$d\eta_{u(t)}=\til{f}_t'(\xi_t)d\xi_t+(\frac\kappa 2-3)\til{f}_t''(\xi_t)dt.$$ 
Since $u'(t)=\til{f}_t'(\xi_t)^2$, so
$$d\eta_s=d\til{\xi}_s+(\frac\kappa 2-3)\til{f}_{v(s)}''(\xi_t)/
\til{f}_{v(s)}'(\xi_t)^2ds,$$ where $\til{\xi}_s=\til{B}(\kappa s)$, $0\le s<p$, and $\til{B}(s)$ is 
another standard Brownian motion on $\R$ started from $0$. Note that $\eta_0=0$. If $\kappa=6$, then $\eta_s=\til{\xi}_s=\til{B}(\kappa s)$, $0\le s<p$. 
Thus $(W(K_{v(s)}))_{0\le s<p}$ has the same 
law as the standard annulus SLE$_{\kappa=6}$ of modulus $p$. So $(K_{v(s)})_{0\le s<p}$
has the same law as $(L_s)_{0\le s<p}$. 

If $\kappa\not=6$, then $d\eta_s=d\til{\xi}_s$ $+$ drift term. 
The remaining part follows from Girsanov's Theorem (\cite{12}). $\Box$

\vskip 2mm
{\bf Remark.} This equivalence implies the a.s.\ existence of annulus SLE trace. 
Suppose $(K_t)$ is an annulus SLE$_\kappa(D;P\to B_2)$. 
If $\kappa\le 4$, the trace
$\beta$ is a simple curve in $D$ such that every $K_t=\beta(0,t]$. If $\kappa>4$ and $B_1$ is 
locally connected, then $\beta$ is a non-simple curve in $D\cup B_1$ such that for every $t$, $D\sem K_t$ 
is the connected component of $D\sem\beta(0,t]$ that has $B_2$ as a boundary component. 

\section{Annulus SLE$_2$ and LERW}
\subsection{Observables for SLE$_2$}
Suppose $D$ is a doubly connected domain of modulus $p$ with boundary components $B_1$ and 
$B_2$, $P$ is a prime end on $B_1$. Let $(K_t)$ be an annulus 
SLE$_2(D;P\to B_2)$ and $\beta$ the corresponding trace. 
Let $D_t=D\sem K_t$, $0\le t<p$. Then $\beta(t,t+\eps)$ determines a prime end
in $D_t$, denoted by $\beta(t_+)$.  
Now consider a positive harmonic function $H_t$ in $D_t$, which has a harmonic conjugate
and satisfies the following properties. As $z\in D_t$ and $z\to B_2$, we have $H_t(z)\to 1$; 
for any neighborhood $V$ of $\beta(t_+)$, as $z\in D_t\sem V$ and $z\to B_1\cup K_t$, 
we have $H_t(z)\to 0$.
The existence of the harmonic conjugate implies that for any smooth Jordan curve, say $\gamma$, that disconnects the two boundary components of $D_t$, we have $\int_\gamma\partial_{\bf n}H_t
ds=0$, where $\bf n$ are normal vectors on $\gamma$ pointed towards $B_1$. 
Now we introduce another positive harmonic function $P_t$
in $D_t$ which satisfies that for any neighborhood $V$ of $\beta(t_+)$, 
as $z\in D_t\sem V$ and $z\to \pa D_t$, we have $P_t(z)\to 0$, and
 $\int_\gamma\partial_{\bf n} P_tds=2\pi$ for any smooth Jordan curve 
$\gamma$ that disconnects the two boundary components of $D_t$. 

\begin{Proposition} For any fixed $z\in D$, $H_t(z)$ and $P_t(z)$, $0\le t<p$, are local 
martingales. \end{Proposition}
{\bf Proof.} By conformal invariance, we may assume that $D=\A_p$, $B_1=\CC_0$, $B_2=\CC_p$,
and $P=1$. So $(K_t,0\le t<p)$ is the standard annulus SLE$_2$ of modulus $p$. 
Let $\chi_t$ and $\vphi_t$, $0\le t<p$,
be the corresponding driving function and conformal maps. Then $\chi_t=\exp(i\xi(t))$
and $\xi(t)=B(2t)$. 
Since $\vphi_t$ maps $D_t$ conformally onto $\A_{p-t}$ and by Proposition 2.1, 
$\vphi_t(\beta(t_+))=\chi_t$, we have 
$$H_t(z)=\Ree \SA_{p-t}(\vphi_t(z)/\chi_t)\mbox{, }\mbox{ and }\mbox{ }P_t(z)=\ln|\vphi_t(z)|+(p-t)H_t(z).$$
We want to use the It\^o's formula. To simplify the computation,
we lift the maps to the covering space. Let $\til{D}_t$, $\til{\A}_r$ and
$\til{\CC}_r$ be the preimages of $D_t$, $\A_r$ and $\CC_r$, respectively, under the map
$z\mapsto e^{iz}$. We may lift $\vphi_t$ to a conformal map $\til{\vphi}_t$ from $\til{D}
_t$ onto $\til{\A}_{p-t}$ so that $\exp(i\til{\vphi}_t(z))=\vphi_t(e^{iz})$, 
$\til{\vphi}_0(z)=z$, and $\til{\vphi}_t(z)$ is continuous in $t$. Let $\til{\SA}_r(z)
=\frac 1i\SA_r(e^{iz})$. Then we have 
$$\partial_t\til{\vphi_t}(z)=\til{\SA}_{p-t}(\til{\vphi_t}(z)-\xi(t)).$$
It is clear that $\til{\SA}_r$ has period $2\pi$, is meromorphic in $\C$ with poles
$\{2k\pi+i2mr:k,m\in\Z\}$, $\Imm\til{\SA}_r\equiv 0$ on $\R\sem\{\mbox{poles}\}$,
and $\Imm\til{\SA}_r\equiv -1$ on $\til{\CC}_r$. It is also easy to check that 
$\til{\SA}_r$ is an odd function, and the principal part of $\til{\SA}_r$ at $0$ is $2/z$. 
So $\til{\SA}_r(z)=2/z+az+O(z^3)$ near $0$, for some $a\in\R$. It is possible to explicit 
this kernel using classical functions in \cite{a}:
$$\til{\SA}_r(z)=2\zeta(z)-\frac{2}\pi\zeta(\pi)z=\frac{1}\pi\frac{\pa_v\theta}\theta
(\frac{z}{2\pi},\frac{ir}{\pi}),$$
where $\zeta$ is the Weierstrass zeta function with basic periods $(2\pi,i2r)$, and
$\theta=\theta(v,\tau)$ is Jacobi's theta function. The following lemma is a direct consequence
of the heat-type differential equation satisfied by $\theta$: $(\pa_v^2-4i\pi\pa_\tau)\theta
=0$. But we prefer a proof using only basic complex analysis. 
The symbols $'$ and $''$ in the lemma denote the first and second derivatives w.r.t.\ $z$.
 
\begin{Lemma} $\partial_r\til{\SA}_r-\til{\SA}_r\til{\SA}_r'-\til{\SA}_r''\equiv 0$. \end{Lemma}
{\bf Proof.} Let $J=\partial_r\til{\SA}_r-\til{\SA}_r\til{\SA}_r'-\til{\SA}_r''$. Then $J$ is odd, has period $2\pi$, takes real values
on $\R\sem\{2k\pi:k\in\Z\}$, and is analytic on $\C\sem\{2k\pi+i2mr:k,m\in\Z\}$. Since near $0$, 
$\til{\SA}_r(z)=2/z+az+O(z^3)$, so $\til{\SA}_r'(z)=-2/z^2+a+O(z^2)$, and $\til{\SA}_r''(z)=4/z^3+O(z)$. Thus $\til{\SA}_r(z)\til{\SA}_r'(z)
+\til{\SA}_r''(z)=O(z)$ near $0$, i.e. $0$ is a removable pole of $\til{\SA}_r\til{\SA}_r'+\til{\SA}_r''$. 
Since $\til{\SA}_r(z)-\frac{1}{i}
\frac{1+e^{iz}}{1-e^{iz}}$ is analytic in a neighborhood of $0$, and 
$\frac{1+e^{iz}}{1-e^{iz}}$ is constant in $t$, so
$0$ is also a removable pole of 
$\partial_r \til{\SA}_r$. Thus $J$ extends analytically at $0$. As $J$ has period $2\pi$, $J$ extends
analytically at $2k\pi$, for all $k\in\Z$. So $J$ is analytic in $\{|\mbox{Im}\,z|<2r\}$. 
The fact that $\Imm \til{\SA}_r\equiv 0$ on $\R\sem\{\mbox{poles}\}$ implies $\Imm J\equiv 0$ 
on $\R$. 

Since $\Imm \til{\SA}_r\equiv -1$ on $\til{\CC}_r=ir+\R$,  we have $\Imm \til{\SA_r}''=\partial_x^2\Imm \til{\SA}_r=
\partial_x\Imm \til{\SA}_r\equiv 0$, and 
$\partial_r\Imm \til{\SA}_r=-\partial_y\Imm \til{\SA}_r=-\partial_x\Ree \til{\SA}_r$ on $\til{\CC}_r$. 
Therefore
$$\Imm(\til{\SA}_r\til{\SA}_r')=\Ree \til{\SA}_r\partial_x\Imm \til{\SA}_r+\Imm \til{\SA}_r\partial_x\Ree \til{\SA}_r=-\partial_x\Ree \til{\SA}_r$$
on $\til{\CC}_r$.
Thus $\mbox{Im}\,J=\Imm\partial_r\til{\SA}_r-\Imm(\til{\SA}_r\til{\SA}_r')-\Imm\til{\SA}_r''\equiv 0$ on 
$\til{\CC}_r$. Now $\Imm J\equiv 0$ on both $\R$ and $ir+\R$, so it has to be zero 
everywhere. It then
follows that $J\equiv C$ for some $C\in\R$. Since $J$ is odd, $C=0$ and $J\equiv 0$. $\Box$
\vskip 2mm
Now we may express $H_t$ and $P_t$ by 
$$H_t(e^{iz})=\Imm\til{\SA}_{p-t}(\til{\vphi_t}(z)-\xi(t))\mbox{, }\mbox{ and }\mbox{ } 
P_t(e^{iz})=\Imm\til{\vphi}_t(z)+(p-t)H_t(z).$$ 
So it suffices to prove that for any $z\in\til{\A}_p$, 
$$M_1(t)=\til{\SA}_{p-t}(\til{\vphi_t}(z)-\xi(t))\mbox{, }\mbox{ and }\mbox{ }
M_2(t)=\til{\vphi}_t(z)+(p-t)M_1(t),$$
$0\le t<p$, are martingales. Using It\^o's formula, we have
$$dM_1(t)=-\partial_r\til{\SA}_{p-t}dt+\til{\SA}_{p-t}'\cdot[d\til{\vphi_t}(z)-d\xi(t)]
+\til{\SA}_{p-t}''dt,$$
where $\partial_r\til{\SA}_{p-t}$, $\til{\SA}_{p-t}'$ and $\til{\SA}_{p-t}''$ are all valued at
$\til{\vphi_t}(z)-\xi(t)$. The last term is the drift term. Note that we use $\kappa=2$ 
here. Since $d\til{\vphi}_t(z)=\til{\SA}_{p-t}(\til{\vphi}_t(z)-\xi(t))dt$, we have
$$dM_1(t)=(-\partial_r\til{\SA}_{p-t}+\til{\SA}_{p-t}'\til{\SA}_{p-t}+\til{\SA}_{p-t}'')dt
-\til{\SA}_{p-t}'d\xi(t)=-\til{\SA}_{p-t}'d\xi(t)$$
by Lemma 3.1. Thus $(M_1(t), 0\le t<p)$ is a local martingale. Now
$$dM_2(t)=\til{\SA}_{p-t}(\til{\vphi_t}(z)-\xi(t))dt+(p-t)dM_1(t)-M_1(t)dt=(p-t)dM_1(t).$$
Thus $(M_2(t),0\le t<p)$ is also a local martingale. $\Box$
\vskip 2mm
{\bf Remark.}
Similar observables also exist for radial and chordal SLE$_2$. For example, let $K_t$ be radial
SLE$_2$ in a simply connected domain $\Omega$, let $H_t$ be the positive harmonic 
function in $\Omega\sem K_t$ which tends to $0$ on $\partial(\Omega\sem K_t)$ except at the ``tip" point
of $K_t$, and normalized so that the value of $H_t$ at the target point is constant $1$.
Then for any fixed $z\in D$, $H_t(z)$, $0\le t<\infty$, is a martingale. This observable 
was mentioned implicitly in the proof of Proposition 3.4 in \cite{6}. As we want to 
define SLE for general domains, we conjecture that such kinds of observables always 
exist for SLE$_2$. 

\subsection{Observables for LERW} 
Let $G=(V,E)$ be a finite or infinite simple connected graph such that $\de(v)<\infty$ for each
$v\in V$. For a function $f$ on $V$, and $v\in V$, let $\Delta_G f(v)=\sum_{w\sim v}(f(w)-f(v))$,
where $w\sim v$ means that $w$ and $v$ are adjacent. A subset $K$ of $V$ is called reachable, if for any 
$v\in V\sem K$, a symmetric 
random walk on $G$ started from $v$ will hit $K$ in finite steps almost 
surely. For subsets $S_1$, $S_2$ and $S_3$ of $V$, let
$\Gamma^{S_3}_{S_1,S_2}$ denote the set of all lattice paths $\gamma=(\gamma_0,\dots,\gamma_n)$ 
such that $\gamma_0\in S_1$, $\gamma_n\in S_2$ and $\gamma_s\in S_3$ for $0<s<n$.
For a finite lattice path $\gamma=(\gamma_0,\dots,\gamma_n)$, write
$$P(\gamma)=1/\prod_{j=0}^n\de(\gamma_j)\mbox{, }\mbox{ }P_0(\gamma)=1/\prod_{j=0}^{n-1}\de(\gamma_j)\mbox{, }
\mbox{ and }\mbox{ }P_1(\gamma)=1/\prod_{j=1}^{n-1}\de(\gamma_j).$$ 
Let $R(\gamma)=(\gamma_n,\dots,\gamma_0)$ be the reversal of $\gamma$, then
$P(R(\gamma))=P(\gamma)$ and $P_1(R(\gamma))=P_1(\gamma)$. If $S_1$, $S_2$ and $S_3$ partition
$V$, $v\in S_3$, then the probability that a random walk on $G$ started from $v$ hits $S_2$ 
before $S_1$ is 
equal to the summation of $P_0(\gamma)$, where $\gamma$ runs over $\Gamma_{v,S_2}^{S_3}$. 

\begin{Lemma} Suppose $A$ and $B$ are disjoint subsets of $V$, and $A\cup B$ is reachable.
Let $f(v)$ be the probability that the random walk on $G$ started from $v$ hits $A$ before $B$. 
Then $f$ is the unique bounded function on $V$ that satisfies $f\equiv 1$ on $A$, $f\equiv 0$ on $B$, and 
$\Delta_G f\equiv 0$ on $C=V\sem(A\cup B)$. Moreover
$\sum_{v\in B}\Delta_Gf(v)=-\sum_{v\in A} \Delta_Gf(v)>0$. \end{Lemma}
{\bf Proof.} The proof is elementary. For the last statement, note that $\sum_{v\in B}\Delta_Gf(v)
=\sum P_1(\gamma)$ where $\gamma$ runs over the non-empty set $\Gamma_{B,A}^{C}$; and 
$-\sum_{v\in A}\Delta_Gf(v)=\sum P_1(\gamma)$ where $\gamma$ runs over $\Gamma_{A,B}^C$. 
The values of the two summations are equal because the reverse map $R$ is a one-to-one correspondence
between $\Gamma_{B,A}^{C}$ and $\Gamma_{A,B}^C$, and $P_1(\gamma)=P_1(R(\gamma))$. $\Box$
\vskip 2mm
Let $L(A,B)=\sum_{v\in B}\Delta_Gf(v)$ for the $f$ in Lemma 3.2. Then $L(A,B)=L(B,A)>0$. 
If any of $A$ or $B$ is a finite set, then we have $L(A,B)<\infty$. 

\begin{Lemma} Let $A$, $B$, $C$ and $f$  be as in Lemma 3.2. Fix $x\in C$. Let $h(v)$ be equal to the 
probability that a simple random walk on $G$ started from $v$ hits $x$ before $A\cup B$. 
Then $$\sum_{v\in A}\Delta_Gh(v)=f(x)(-\Delta_Gh(x)).$$  \end{Lemma}
{\bf Proof.} From the proof of Lemma 3.2, we have 
$$f(x)=\sum_{\alpha\in\Gamma^C_{x,A}}P_0(\alpha)=\sum_{\beta\in\Gamma_{x,x}^C}P(\beta)
\sum_{\gamma\in\Gamma^{C\sem\{x\}}_{x,A}}P_1(\gamma)=\sum_{\beta\in\Gamma_{x,x}^C}P(\beta)
\sum_{v\in A}\Delta_Gh(v),$$
and 
$$1=\sum_{\alpha\in\Gamma^C_{x,A\cup B}}P_0(\alpha)=\sum_{\beta\in\Gamma_{x,x}^C}P(\beta)
\sum_{\gamma\in\Gamma^{C\sem\{x\}}_{x,A\cup B}}P_1(\gamma)=\sum_{\beta\in\Gamma_{x,x}^C}P(\beta)
(-\Delta_G h(x)).$$
So we proved this lemma. $\Box$

\begin{Lemma} Let $A$, $B$, $C$ and $f$  be as in Lemma 3.2. Suppose $L(A,B)<\infty$. 
Fix $x\in C$ such that $f(x)>0$. Then there is a unique bounded function $g$ on $V$ such that
$g\equiv 1$ on $A$; $g\equiv 0$ on $B$; $\Delta_Gg\equiv 0$ on $C\sem\{x\}$; and 
$\sum_{v\in A}\Delta_Gg(v)=0$. Moreover, such $g$ is non-negative and satisfies
$\sum_{v\in B\cup\{x\}}\Delta_Gg(v)=0$ and $\Delta_Gg(x)=-L(A,B)/f(x)$. \end{Lemma}
{\bf Proof}. Suppose $g$ satisfies the first group of properties. Let $I=g-f$. Then $I$ is bounded, 
$I\equiv 0$ on $A\cup B$ and $\Delta_GI\equiv 0$ on $C\sem\{x\}$. Thus $I(v)=I(x)h(v)$, where $h$ 
is as in Lemma 3.3. Then by Lemma 3.2 and 3.3, 
$$0=\sum_{v\in A}\Delta_Gg(v)=\sum_{v\in A}\Delta_G(I+f)(v)=
-I(x)f(x)\Delta_Gh(x)-L(A,B).$$
Thus $I(x)=L(A,B)/(-f(x)\Delta_Gh(x))$ is uniquely determined. Therefore $g$ is unique. 

On the other hand, if we define $g=f+h L(A,B)/(-f(x)\Delta_Gh(x))$, then from the last paragraph, we see 
that $g$ satisfies
the first group of properties. Since $f$ and $h$ are non-negative, and 
$-\Delta_Gh(x)=L(x,A\cup B)>0$ by
Lemma 3.2, so $g$ is also non-negative. By Lemma 3.2 and 3.3,
$$\sum_{v\in B\cup\{x\}}\Delta_Gg(v)=L(A,B)+\Delta_Gf(x)+
\sum_{v\in B\cup\{x\}}\Delta_Gh(v)L(A,B)/(-f(x)\Delta_Gh(x))$$
$$=L(A,B)-\sum_{v\in A}\Delta_Gh(v)L(A,B)/(-f(x)\Delta_Gh(x))=L(A,B)-L(A,B)=0.$$
Finally, $\Delta_Gg(x)=\Delta_Gh(x)\cdot L(A,B)/(-f(x)\Delta_Gh(x))=-L(A,B)/f(x)$. $\Box$
\vskip 2mm
From now on, let $D$ be a doubly connected domain with boundary components $B_1$ and $B_2$, and
satisfies $0\in B_1$ and $(0,a]\subset D$ for some $a>0$. We use the symbols $D^\delta$ and LERW 
defined in Section 1.3. Note that $D^\delta$ may not be connected. To apply
the lemmas in above, we need to modify $D^\delta$ a little bit. Let $\cal P$ denote the set 
of all lattice paths on $D^\delta$ from $\delta$ to some boundary vertex whose vertices
are inside $D$ except the last vertex.  
Every path of $\cal P$ can be viewed as a subgraph of $D^\delta$.
Let $\til{D^\delta}$ be the union of all paths in $\cal P$ as a subgraph of $D^\delta$.
Then $\til{D^\delta}$ is a connected graph. And if we replace $D^\delta$ by $\til{D^\delta}$
in the definition of LERW in Section 1.3, we will get the same 
LERW. So we can consider $\til{D^\delta}$ instead of $D^\delta$. For simplicity of notations,
we write $D^\delta$ for $\til{D^\delta}$. 

By the definition, any two vertices of $D^\delta$ on $\partial D$ are not adjacent,
so the neighbors of boundary vertices of $D^\delta$ are those vertices lie in $D$, are in $\delta\Z^2$
and has exactly $4$ neighbors. It follows that if any $B_j$ is bounded, then there are finitely  many vertices that
lie on $B_j$. On the other hand, $B_1$ and $B_2$ can't be both unbounded. Now we denote
$$E_{-1}^\delta=V(D^\delta)\cap B_1\mbox{, }\mbox{ }F^\delta=V(D^\delta)\cap B_2\mbox{, }\mbox{ and }\mbox{ }
N_{-1}^\delta=V(D^\delta)\cap D.$$ 
Then at least one of $E_{-1}^\delta$ and $F^\delta$ is a finite set. 
Write LERW as $y=(y_0,\dots,y_\up)$, where $y_0=\delta$ and $y_\up\in B_2$. For $0\le j<\up$, let 
$$E_j^\delta=E_{-1}^\delta\cup\{y_0,\dots,y_j\}\mbox{, }\mbox{ and }\mbox{ }
N_j^\delta=N_{-1}^\delta\sem\{y_0,\dots,y_j\}.$$ 
Then $E_j^\delta$, $N_j^\delta$ and $F^\delta$ partition $V(D^\delta)$, for $-1\le j<\up$.
The fact that the lattice $\Z^2$ is recurrent easily implies that $E_j^\delta\cup F^\delta$ 
is reachable in $D^\delta$. 
Since one of $E_{j}^\delta$ and $F^\delta$ is a finite set, we have $L(E_j^\delta,F^\delta)<\infty$ for $-1\le j<\up$.
For $-1\le j<\up$, let $f_j$ be the $f$ in Lemma 3.2 with $G=D^\delta$, $A=F^\delta$ 
and $B=E_j^\delta$. For 
$0\le j<\up$, since $(y_j,\dots,y_\up)$ is a lattice path from $y_j$ to $F^\delta$ not 
passing through $E_{j-1}^\delta$, we have $f_{j-1}(y_j)>0$.  
Let $g_j$ be the $g$ in Lemma 3.4 with $G=D^\delta$, $A=F^\delta$, $B=E_{j-1}^\delta$, 
and $x=y_j$, for $0\le j<\up$.

\begin{Lemma} Conditioned on the event that $y_j=w_j$, $0\le j\le k$, and $k<\up$, the 
probability that $y_{k+1}=u$ is $f_k(u)/\sum_{v\sim w_k}f_k(v)$ if $u\sim w_k$; and 
is zero if $u\not\sim w_k$. \end{Lemma}
{\bf Proof.} This result is well known. See \cite{2} for details. $\Box$

\begin{Proposition}
Let $\lin{F^\delta}$ be the union of $F^\delta$ and the set of vertices of $D^\delta$ that are adjacent to $F^\delta$. 
Fix a vertex $v_0$ of $D^\delta$. Conditioned on the event that $y_j=w_j$, $0\le j\le k$, 
$w_k\not\in\lin{F^\delta}$, and $f_k(v_0)>0$, the expectation of $g_{k+1}(v_0)$ is equal to $g_k(v_0)$, which is
determined by $w_j$, $0\le j\le k$. Thus $g_k(v_0)$ is a discrete martingale up to the first time
$y_k$ hits $\lin{F^\delta}$, or $E_k^\delta=E_{-1}^\delta\cup\{y_0,\dots,y_k\}$ 
disconnects $v_0$ from $F^\delta$ in $D^\delta$. \end{Proposition}
{\bf Proof.} Let $S$ be the set of $v$ such that $v\sim w_k$ and $f_k(v)>0$. By lemma 3.5, 
the conditional probability that $y_{k+1}=u$ is $f_k(u)/\sum_{v\in S}f_k(v)$ for $u\in S$. 
For $v\in S$, let $g_{k+1}^v$ be the $g$ in Lemma 3.4 with $G=D^\delta$, $A=F^\delta$, 
$B=E_k^\delta$ and $x=v$.
Then with probability $f_k(u)/\sum_{v\in S}f_k(v)$, $g_{k+1}=g_{k+1}^u$. Thus the conditional
expectation of $g_{k+1}(v_0)$ is equal to $\til{g}_k(v_0)$, where 
$$\til{g}_k(v):=\sum_{u\in S}f_k(u)g_{k+1}^u(v)/\sum_{u\in S}f_k(u).$$ 
Then $\til{g}_k\equiv 0$ on $E_k^\delta$, 
$\equiv 1$ on $F^\delta$; $\Delta\til{g}_k\equiv 0$ on $N_k^\delta\sem S$, and $\sum_{v\in F^\delta}\Delta\til{g}_k(v)=0$.
And $$\Delta\til{g}_k(v)=\frac{f_k(v)\Delta g_{k+1}^v(v)}{\sum_{u\in S} f_k(u)}
=-\frac{L(E_k^\delta,F^\delta)}{\sum_{u\in S} f_k(u)}\mbox{, }\mbox{ }\mbox{ }\forall v\in S,$$ 
by Lemma 3.4. Now define $\ha{g}_k$ on 
$V(D^\delta)$ such
that $\ha{g}_k(w_k)=L(E_k^\delta,F^\delta)/\sum_{u\in S} f_k(u)$; for those $v\in N_k^\delta$ 
such that $f_k(v)=0$,
define $\ha{g}_k(v)$ to be $\ha{g}_k(w_k)$ times the probability that a simple random walk on $D^\delta$ started from 
$v$ hits $w_k$ before $E_{k-1}^\delta$; and let $\ha{g}_k(v)=\til{g}_k(v)$ for other $v\in V(D^\delta)$. 
Then $\Delta \ha{g}_k\equiv 0$ on $N_k^\delta$, $\ha{g}_k\equiv 0$ on $E_k^\delta\sem\{w_k\}$, 
and $\ha{g}_k
\equiv 1$ on $F^\delta$. Since $w_k\not\in\lin{F^\delta}$, and for $v\in N_k^\delta$ such that
$f_k(v)=0$ we have $v\not\in\lin{F^\delta}$, so 
$\sum_{v\in F^\delta}\Delta\ha{g}_k(v)=\sum_{v\in F^\delta}\Delta
\til{g}_k(v)=0$. Now $\ha{g}_k$ satisfies all properties of $g_k$. The uniqueness of $g_k$ implies that
$\ha{g}_k\equiv g_k$. Since $f_k(v_0)>0$, we have $g_k(v_0)=\ha{g}_k(v_0)=\til{g}_k(v_0)$. $\Box$
\vskip 2mm
{\bf Remark 1}. The observable $g_k$ corresponds to $H_t$ in Proposition 3.1.
We may define another kind of observables $q_k$ to be the bounded function on the vertices of $D^\delta$ such
that $q_k\equiv 0$ on $E_{k-1}\cup F$, $\Delta q_k\equiv 0$ on $N_k$, and $\sum_{v\in F}\Delta
q_k(v)=2\pi=-\sum_{v\in E_k}\Delta q_k(v)$. Then Proposition 3.2 still holds if $g_k$ is 
replaced by $q_k$, and $q_k$ corresponds to $P_t$ in Proposition 3.1. 
The definition of $q_k$ does not need the fact that $L(E_k,F)<\infty$.
We may also use $q_k$ to do the approximation. 

{\bf Remark 2}. Suppose $\alpha$ is a Jordan curve in $D$ which disconnects $E_k^\delta$ 
from $F$ and
does not pass through any vertex of $D^\delta$. Denote $D_j$ the component of $D\sem\alpha$ that has 
$B_j$ as part of boundary, $j=1,2$. We also suppose that $y_0$ through $y_k$ are in $D_1$. Let $S$ 
be the set of vertex pair $(v,w)$ such that $v\in D_1$, $w\in D_2$, and $v\sim w$. 
From the fact that $\Delta_{D^\delta} g_k\equiv 0$ on $V(D^\delta)\cap D$, we conclude
$\sum_{(v,w)\in S}(g_k(v)-g_k(w))=0$. Similarly, $\sum_{(v,w)\in S}(q_k(v)-q_k(w))=2\pi$. 
\vskip 2mm
Now suppose $\alpha_1$ and $\alpha_2$ are two disjoint Jordan curves in $D$ such that $\alpha_j$
disconnects $\alpha_{3-j}$ from $B_j$, $j=1,2$. For $j=1,2$, let $U_j$ be the subdomain of $D$
bounded by $\alpha_j$ and $B_j$, and $V_j^\delta=V(D^\delta)\cap U_j$. Let $L^\delta$ be the set
of simple lattice paths of the form $w=(w_{-1},w_0,\dots,w_k)$, $k\ge 0$ such that $w_{-1}\in B_1$,
$w_0,\dots,w_k\in V_1^\delta$, and there is some lattice path from the last vertex $P(w):=w_k$ to $B_2$ 
without passing $w_0,\dots,w_{k-1}$, and vertices on $B_1$. For $w\in L^\delta$, denote 
$$E^\delta_w=E^\delta_{-1}\cup\{w_0,\dots,w_k\}\mbox{, }\mbox{ and }\mbox{ }
N^\delta_w=N^\delta_{-1}\sem\{w_0,\dots,w_k\}.$$ 
Let $g_w$ be the $g$ in Lemma 3.4 
with $G=D^\delta$, $A=F^\delta$, $B=E^\delta_w\sem\{P(w)\}$, and $x=P(w)$. Now define $D_w=D\sem\cup_{j=0}
^k[w_{j-1},w_j]$. Let $u_w$ be the non-negative harmonic function in $D_w$ 
whose harmonic conjugates exist, and whose continuation is constant $1$ on $B_2$,
and constant $0$ on $\cup_{j=0}^k[w_{j-1},w_j]\cup B_1$ except at $P(w)$. 
The existence of the harmonic conjugates implies that 
$\int_\alpha\partial_{\bf n} u_wds=0$ for any 
smooth Jordan curve $\alpha$ that disconnects $B_2$ from $\cup_{j=0}
^k[w_{j-1},w_j]\cup B_1$.
It is intuitive to guess that $g_w$ should be close to $u_w$.
In fact, we have the following proposition. The proof is postponed to Section 5.

\begin{Proposition} Given any $\eps>0$, there is $\delta(\eps)>0$ such that if $0<\delta<\delta(\eps)$
and $w\in L^\delta$, then $|g_w(v)-u_w(v)|<\eps$, for any $v\in V_2^\delta$. \end{Proposition}

\subsection{Convergence of the driving process}
Fix some small $\delta>0$. We write LERW on $D^\delta$ by $y=(y_0,\dots,y_\up)$ as in Section 3.2. 
Let $y_{-1}=0$. Extend $y$ to be a map from $[-1,\up]$ into $\lin{D}$ such that $y$ is linear on $[j-1,j]$ for each
$0\le j\le\up$. It clear that $y(-1,s]$, $-1\le s<\up$, is a Loewner chain in $D$ on
$B_1$. And $y(-1,s]$ approaches $B_2$ as $s\nearrow\up$. For $-1\le s<\up$, let $T(s)=C_{D,B_2}(y(-1,s])$, 
then $T$ is a continuous increasing function, and maps $[-1,\up)$ onto $[0,p)$, where $p=M(D)$.  
Let $S:[0,p)\to[-1,\up)$ be the inverse of $T$. Let $\beta(t)=y(S(t))$, and $K_t=\beta(0,t]$, for $0\le t<p$. 
Suppose $W$ maps $D$ conformally onto $\A_p$ so that $W(0_+)=1$, i.e., $W(x)\to 1$ as $x\in\R^+$ and $x\to 0$. 
Then $(W(K_t),0\le t<p)$ is a Loewner chain in $\A_p$ on $\CC_0$ such that $C_{A_p,C_p}(W(K_t))=t$. 
By Proposition 2.1, $W(K_t)$, $0\le t<p$, are the standard annulus LE hulls of modulus $p$
driven by some continuous $\chi_t$, $0\le t<p$, on $\CC_0$. 
Let $(\vphi_t,0\le t<p)$ be the corresponding standard annulus LE maps.
Since $W(\beta(t))\to 1$ as $t\to 0$, $\chi_0=1$. We may write $\chi_t=e^{i\xi_t}$, so that $\xi_0=0$, and $\xi_t$ is continuous in $t$. We want to prove that the law of 
$(\xi_t)_{0\le t<p}$, which depends on $\delta$, converges to the law of $(B(2t))_{0\le t<p}$. 

For $a<b$, let $\A_{a,b}$ be the annulus bounded by $\CC_a$ and $\CC_b$. 
For any $0<q<p$, there is a smallest $l(p,q)\in(0,p)$ such that
if $K$ is a hull in $\A_p$ on $\CC_0$ with the capacity (w.r.t.\ 
$\CC_p$) less than $q$, then $K$ does not intersect $\A_{l(p,q),p}$. 
Using the fact that for any
$0<s\le r$, $\Ree\SA_r$ attains its unique maximum and minimum on $\lin{\A_{s,r}}$ at
$e^{-s}$ and $-e^{-s}$, respectively, it is not hard to derive the following Lemma.

\begin{Lemma} Fix $0<q<p$, let $r\in(l(p,q),p)$. There are $\iota\in(0,1/2)$ and $M>0$ 
depending on $p$, $q$ and $r$, which satisfy the following properties. Suppose $\vphi_t$, 
$0\le t<p$, are
some standard annulus LE maps of modulus $p$ driven by $\chi_t$, $0\le t<p$. 
Then we have  
$|\partial_z\SA_{p-t}(\vphi_t(z)/\chi_t)|\le M$, for all $t\in[0,q]$ and $z\in\A_{r,p}$. 
Moreover, $$\A_{\iota(p-t),p-t}\supset\vphi_t(\A_{r,p})\supset
\A_{(1-\iota)(p-t),p-t}\mbox{, }\mbox{ }\forall t\in[0,q].$$
\end{Lemma}

Now fix $q_0\in(0,p)$. Let $q_1=(q_0+p)/2$. Choose $p_1\in(l(p,q_1),p)$, and let $p_2=(p_1+p)/2$. Denote $\alpha_j=W^{-1}(\CC_{p_j})$, $j=1,2$. Then 
$\alpha_1$ and $\alpha_2$ are disjoint Jordan curves in $D$ such that $\alpha_j$
disconnects $\alpha_{3-j}$ from $B_j$, $j=1,2$.
Let $n_\infty=\lceil S(q_0)\rceil$, where $\lceil x\rceil$ is the smallest integer that is not less than $x$. 
Then $n_\infty$ is a stopping time w.r.t.\ $\{{\cal F}_k\}$, where ${\cal F}_k$ 
denotes the $\sigma$-algebra generated by $y_0$, $y_1$, $\dots$, $y_{k\land\up}$.
For $0\le k\le n_\infty-1$, $T(k)\le q_0<q_1$, so  
from the choice
of $p_1$, we see that $W(y_k)$ lies in the domain bounded by $\CC_{p_1}$ and $\CC_0$, so $y_k$ 
lies
in the domain bounded by $B_1$ and $\alpha_1$. Note that $y_{-1}=0\in B_1$. 
So for $-1\le k\le n_\infty-1$, 
if $\delta$ is small, then $[y_k,y_{k+1}]$ can be disconnected from $B_2$ by an annulus centered
at $y_k$ with inner radius $\delta$ and outer radius $dist(\alpha_1, B_2)$. So as $\delta\to0$, the conjugate extremal distance between $B_2$ and $[y_k,y_{k+1}]$ in 
$D_{y^k}=D\sem\cup_{0\le j\le k}[y_{j-1},y_j]$ (the extremal length of the family of 
rectifiable curves in $D_{y^k}$ that disconnect $B_2$ from  $[y_k,y_{k+1}]$,
see \cite{1}) tends to $0$, uniformly in $-1\le k\le n_\infty-1$. It then follows that 
$T(k+1)-T(k)$ and $\max\{|\xi_t-\xi_{T(k)}|:T(k)\le t\le T(k+1)\}$ tend to $0$ as $\delta\to 0$, uniformly
in $-1\le k\le n_\infty-1$.
Since $T(n_\infty-1)\le q_0$, we may choose $\delta$ small enough such that $T(n_\infty)<q_1$. 
We now use the symbols in the last part of Section 3.2 for Jordan
curves $\alpha_1$ and $\alpha_2$ defined here. For $0\le k\le n_\infty$, 
let $y^k=(y_{-1},y_0,\dots,y_k)\in L^\delta$.
Then $g_{y^k}=g_k$. By Proposition 3.2, for any fixed $v\in V_2^\delta$, 
$g_k(v)$, $0\le k\le n_\infty$, is a discrete martingale w.r.t.\ $\{{\cal F}_k\}$.

Now fix $d>0$. Define a non-decreasing sequence $(n_j)_{j\ge 0}$ inductively. Let $n_0=0$. Let $n_{j+1}$ be the 
first integer $n\ge n_j$ such that $T(n)-T(n_j)\ge d^2$, or $|\xi_{T(n)}-\xi_{T(n_j)}|\ge d$, or $n\ge n_\infty$, 
whichever comes first. Then $n_j$'s are stopping times w.r.t.\ $\{{\cal F}_k\}$, and they are bounded above by
$n_\infty$. If we let $\delta$ be smaller than some constant depending on $d$, 
then $T(n_{j+1})-T(n_j)\le 2d^2$ and $|\xi_{T(s)}-\xi_{T(n_j)}|\le 2d$ for all $s\in[n_j,n_{j+1}]$ and $j\ge0$. 
Let ${\cal F}_j'={\cal F}_{n_j}$. Then for any $v\in V_2^\delta$, $\{g_{n_j}(v):0\le j<\infty\}$ is a discrete martingale w.r.t.\ $\{{\cal F}_j'\}$. 
Since $\vphi_{T(k)}\circ W$ maps $D_{y^k}$ conformally onto $\A_{p-T(k)}$ and takes
$y_k=P(y^k)$ to $\chi_{T(k)}$, we have
$$u_{y^k}(z)=\Ree\SA_{p-T(k)}(\vphi_{T(k)}\circ W(z)/\chi_{T(k)}).$$ 
By Proposition 3.3, for any $z\in W(V^\delta_2)$ and $0\le j\le k$, \begin{equation*}
\EE[\Ree\SA_{p-T(n_k)}(\vphi_{T(n_k)}(z)/\chi_{T(n_k)})|
{\cal F}_j']=\Ree\SA_{p-T(n_j)}(\vphi_{T(n_j)}(z)/\chi_{T(n_j)})+o_\delta(1).\end{equation*}
As $\delta$ tends to $0$, the set $W(V^\delta_2)$ tends to be dense in $\A_{p_2,p}$. So
for any $z\in\A_{p_2,p}$, there is some $z_0\in W(V^\delta_2)$ such that $|z-z_0|=o_\delta(1)$.
Note that $T(n_j)\le T(n_k)\le T(n_\infty)\le q_1$ for $0\le j\le k$. 
Using the boundedness of the derivative in Lemma 3.6 with $q=q_1$ and $r=p_2$, 
we then have that for all $z\in\A_{p_2,p}$, 
$$\EE[\mbox{Re}\,\SA_{p-T(n_k)}(\vphi_{T(n_k)}(z)/\chi_{T(n_k)})|
{\cal F}_j']=\mbox{Re}\,\SA_{p-T(n_j)}(\vphi_{T(n_j)}(z)/\chi_{T(n_j)})+o_\delta(1).$$

Now consider the maps in the covering space. 
We use the notations in Section 3.1. And let $\til{\A}_{a,b}$ be the preimage of 
$\A_{a,b}$ under the map $z\mapsto e^{iz}$. Then we have
\begin{equation} \EE[\mbox{Im}\,\til{\SA}_{p-T(n_k)}(\til{\vphi}_{T(n_k)}(z)-\xi_{T(n_k)})|
{\cal F}_j']=\mbox{Im}\,\til{\SA}_{p-T(n_j)}(\til{\vphi}_{T(n_j)}(z)-\xi_{T(n_j)})+o_\delta(1).\end{equation}
In Lemma 3.6, let $q=q_1$ and $r=p_2$, then we have some $\iota\in(0,1/2)$ such that  
\begin{equation}\til{\A}_{\iota(p-t),p-t}
\supset\til{\vphi}_t(\til{\A}_{p_2,p})\supset\til{\A}_{(1-\iota)(p-t),p-t},\end{equation} 
for $0\le t\le q_1$.

\begin{Proposition} There are an absolute constant $C>0$ and a constant $\delta(d)>0$ 
such that if $\delta<\delta(d)$, then for all $j\ge 0$,
$$|\EE[\xi_{T(n_{j+1})}-\xi_{T(n_j)}|{\cal F}_j']|\le Cd^3\mbox{, and}$$
$$|\EE[(\xi_{T(n_{j+1})}-\xi_{T(n_j)})^2/2-(T(n_{j+1})-T(n_j))|{\cal F}_j']|\le Cd^3.$$
\end{Proposition}
{\bf Proof}. Fix some $j\ge 0$. Let $a=T(n_j)$ and $b=T(n_{j+1})$. 
Then $0\le a\le b\le q_1$. And if $\delta$ is less than some $\delta_1(d)$,
we have $|b-a|\le 2d^2$ and $|\xi_c-\xi_a|\le 2d$, for any $c\in[a,b]$.  
Now suppose $z\in{\til{\A}_{p_2,p}}$, and consider 
$$I:=\til{\SA}_{p-b}(\til{\vphi}_b(z)-\xi_b)-\til{\SA}_{p-a}(\til{\vphi}_a(z)-\xi_a).$$
Then $I=I_1+I_2$, where
$$I_1:=\til{\SA}_{p-b}(\til{\vphi}_b(z)-\xi_b)-\til{\SA}_{p-b}(\til{\vphi}_a(z)-\xi_a),$$ 
$$I_2:=\til{\SA}_{p-b}(\til{\vphi}_a(z)-\xi_a)-\til{\SA}_{p-a}(\til{\vphi}_a(z)-\xi_a).$$
Then for some $c_1\in[a,b]$, $I_1=I_3+I_4+I_5$, where
$$I_3:=\til{\SA}_{p-b}'(\til{\vphi}_a(z)-\xi_a)[(\til{\vphi}_b(z)-\til{\vphi}_a(z))-(\xi_b-\xi_a)],$$
$$I_4:=\til{\SA}_{p-b}''(\til{\vphi}_a(z)-\xi_a)[(\til{\vphi}_b(z)-\til{\vphi}_a(z))-(\xi_b-\xi_a)]^2/2,$$
$$I_5:=\til{\SA}_{p-b}'''(\til{\vphi}_{c_1}(z)-\xi_{c_1})[(\til{\vphi}_b(z)-\til{\vphi}_a(z))-(\xi_b-\xi_a)]^3/6.$$
And for some $c_2\in[a,b]$, we have
\begin{equation}I_2=-\partial_r\til{\SA}_{p-b}(\til{\vphi}_{a}(z)-\xi_a)(b-a)
+\partial_r^2\til{\SA}_{p-c_2}(\til{\vphi}_a(z)-\xi_a)(b-a)^2/2. \end{equation}
Now for some $c_3\in[a,b]$, we have\begin{equation}
\til{\vphi}_b(z)-\til{\vphi}_a(z)=\partial_r\til{\vphi}_{c_3}(z)(b-a)
=\til{\SA}_{p-c_3}(\til{\vphi}_{c_3}(z)-\xi_{c_3})(b-a).\end{equation}
For some $c_4\in[c_3,b]$, we have\begin{equation}
\til{\SA}_{p-c_3}(\til{\vphi}_{c_3}(z)-\xi_{c_3})=\til{\SA}_{p-b}(\til{\vphi}_{c_3}(z)-\xi_{c_3})+
\partial_r\til{\SA}_{p-c_4}(\til{\vphi}_{c_3}(z)-\xi_{c_3})(b-c_3).\end{equation}
For some $c_5\in[a,c_3]$, we have
$$\til{\SA}_{p-b}(\til{\vphi}_{c_3}(z)-\xi_{c_3})=\til{\SA}_{p-b}(\til{\vphi}_{a}(z)-\xi_{a})$$
\begin{equation}+\til{\SA}_{p-b}'(\til{\vphi}_{c_5}(z)-\xi_{c_5})
[(\til{\vphi}_{c_3}(z)-\til{\vphi}_a(z))-(\xi_{c_3}-\xi_a)].\end{equation}
Once again, there is $c_6\in[a,c_3]$ such that \begin{equation}
\til{\vphi}_{c_3}(z)-\til{\vphi}_a(z)=\partial_r\til{\vphi}_{c_6}(z)(c_3-a)=
\til{\SA}_{p-c_6}(\til{\vphi}_{c_6}(z)-\xi_{c_6})(c_3-a).\end{equation}

We have the freedom to choose $d$ arbitrarily small. Now suppose $d<(1-\iota)(p-q_1)/2$.
Then $$p-a\le p-b+2d\le (p-b)+(1-\iota)(p-q_1)\le (2-\iota)(p-b).$$
Thus for any $m\le M\in[a,b]$, $p-m\le(2-\iota)(p-M)$. By formula (3.2),
$$\til{\vphi}_{m}(z)-\xi_m\in{\til{\A}_{\iota(p-m),p-m}}\subset{\til{\A}_{\iota(p-M),
(2-\iota)(p-M)}}.$$
So the values of $\til{\SA}_{p-M}$, $\partial_r\til{\SA}_{p-M}$, $\partial_r^2\til{\SA}_{p-M}$, $\til{\SA}_{p-M}'$, $\til{\SA}_{p-M}''$ and
$\til{\SA}_{p-M}'''$ at $\til{\vphi}_{m}(z)-\xi_{m}$ are uniformly bounded. 
In formula (3.3), consider $m=a$ and $M=c_2$. Since $|b-a|\le 2d^2$, we have
$$I_2=-\partial_r\til{\SA}_{p-b}(\til{\vphi}_a(z)-\xi_a)(b-a)+O(d^4).$$
Similarly, formula (3.7) implies  
$$\til{\vphi}_{c_3}(z)-\til{\vphi}_a(z)=O(c_3-a)=O(d^2).$$ 
This together with formulae (3.5),(3.6) and $\xi_{c_3}-\xi_a=O(d)$ implies that
$$\til{\SA}_{p-c_3}(\til{\vphi}_{c_3}(z)-\xi_{c_3})=\til{\SA}_{p-b}(\til{\vphi}_{a}(z)-\xi_{a})+O(d).$$ 
By formula (3.4), we have $$\til{\vphi}_b(z)-\til{\vphi}_a(z)=\til{\SA}_{p-b}(\til{\vphi}_{a}(z)-\xi_{a})(b-a)+O(d^3)=O(d^2).$$
Thus  $I_5=O(d^3)$, 
$$I_4=\til{\SA}_{p-b}''(\til{\vphi}_a(z)-\xi_a)(\xi_b-\xi_a)^2/2+O(d^3)\mbox{, }\mbox{ and}$$ 
$$I_3=\til{\SA}'_{p-b}(\til{\vphi}_a(z)-\xi_a)[\til{\SA}_{p-b}
(\til{\vphi}_a(z)-\xi_a)(b-a)-(\xi_b-\xi_a)]+O(d^3).$$
Note that $I=I_2+I_3+I_4+I_5$. Using Lemma 3.1, we get 
$$I=\til{\SA}_{p-b}''(\til{\vphi}_a(z)-\xi_a)[(\xi_b-\xi_a)^2/2-(b-a)]$$
$$-\til{\SA}'_{p-b}(\til{\vphi}_a(z)-\xi_a)(\xi_b-\xi_a)+O(d^3).$$
By formula (3.1), if $\delta$ is smaller than some $\delta_2(d)$, then the conditional expectation of 
$$\mbox{Im}\,\til{\SA}_{p-b}''(\til{\vphi}_a(z)-\xi_a)[(\xi_b-\xi_a)^2/2-(b-a)]
-\mbox{Im}\,\til{\SA}'_{p-b}(\til{\vphi}_a(z)-\xi_a)[\xi_b-\xi_a]$$
w.r.t.\ ${\cal F}'_j$ is bounded by $C_1d^3$. 

By formula (3.2), for any $w\in{\til{\A}_{(1-\iota)(p-a),p-a}}$, the conditional expectation of\begin{equation}
\mbox{Im}\,\til{\SA}_{p-b}''(w)[(\xi_b-\xi_a)^2/2-(b-a)]-\mbox{Im}\,\til{\SA}'_{p-b}(w)[\xi_b-\xi_a]\end{equation}
w.r.t ${\cal F}'_j$ is bounded by $C_1d^3$, if $\delta$ is small enough (depending on $d$).

Now suppose $d<(p-q_1)\iota/(4-4\iota)$. Then $$(1-\iota)(p-a)<(1-\iota/2)(p-b)<p-a.$$ 
Thus $i(1-\iota/2)(p-b)\in\til{\A}_{(1-\iota)(p-a),p-a}$. We may check
$$\mbox{Im}\,\til{\SA}_{p-b}''(i(1-\iota/2)(p-b))>0\mbox{, }\mbox{ and }\mbox{ }
\Imm\til{\SA}_{p-b}'(i(1-\iota/2)(p-b))=0.$$
So we can find $C_2>0$ such that for all $b\in[0,q_1]$,
$\mbox{Im}\,\til{\SA}_{p-b}''(i(1-\iota/2)(p-b))>C_2$. 
Let $w=i(1-\iota/2)(p-b)$ in formula (3.8), then we get
$$|\EE[(\xi_b-\xi_a)^2/2-(b-a)|{\cal F}_j']|\le C_3d^3.$$
Since $\mbox{Im}\,\til{\SA}_{p-b}''(w)$ is uniformly bounded on $\til{\CC}_{(1-\iota/2)(p-b)}$, so
for all $w\in\til{\CC}_{(1-\iota/2)(p-b)}$,
\begin{equation}\mbox{Im}\,\til{\SA}_{p-b}'(w)|\EE[\xi_b-\xi_a|{\cal F}_j']|\le C_4d^3.\end{equation}
We may check that 
$$x_b:=\mbox{Im}\,\til{\SA}_{p-b}(\pi+i(1-\iota/2)(p-b))-\mbox{Im}\,\til{\SA}_{p-b}(i(1-\iota/2)(p-b))>0.$$ So $x_b$ is greater than some absolute constant $C_5>0$ for $b\in[0,q_1]$.
Then there exists $w_b\in\til{\CC}_{(1-\iota/2)(p-b)}$ such that 
$$|\mbox{Im}\,\til{\SA}_{p-b}'(w_b)|=|\partial_x\mbox{Im}\,\til{\SA}_{p-b}(w_b)|=x_b/\pi
\ge C_5/\pi.$$
Plugging $w=w_b$ in formula (3.9), we then have $|\EE[\xi_b-\xi_a|{\cal F}_j']|\le C_6d^3$. $\Box$

\vskip 2mm
The following Theorem about the convergence of the driving process can be deduced from 
Proposition 3.4 
by using the Skorokhod Embedding Theorem. It is very similar to Theorem 3.6 in \cite{6}. 
So we omit the proof.

\begin{Theorem} For every $q_0\in(0,p)$ and $\eps>0$ there is a $\delta_0>0$ depending on $q_0$ and $\eps$ such that for $\delta<\delta_0$ there is a coupling of the processes $\xi_t$ 
and $B(2t)$ such that 
$$\mbox{\bf P}[\sup\{|\xi_t-B(2t)\}|:t\in [0,q_0]\}>\eps]<\eps.$$\end{Theorem}

\subsection{Convergence of the trace}
In this subsection, we will prove Theorem 1.2. We use symbols $y^\delta$, $\beta^\delta$, 
$K^\delta$ and $\chi^\delta$ to emphasize the fact that they depend on $\delta$. Let 
$(K^0_t,0\le t<p)$
be the annulus SLE$_2$ in $D$ from $0_+$ to $B_2$. Let $\beta^0:(0,p)\to D$ be the
corresponding trace.

First, we need two well-known lemmas about simple random walks on $\delta\Z^2$.
We use the superscript $\#$ to denote the spherical metric.

\begin{Lemma} Suppose $v\in\delta\Z^2$ and $K$ is a connected set on the plane that has Euclidean 
(spherical, resp.) diameter at least $R$. Then the probability that a simple random walk on $\delta\Z^2$ started
 from $v$ will exit $\B(v;R)$ ($\B^\#(v;R)$, resp.) before using an edge of 
$\delta\Z^2$ that intersects $K$ is at most $C_0((\delta+dist(v,K))/R)^{C_1}$ 
($C_0((\delta+dist^\#(v,K))/R)^{C_1}$, resp.) for some absolute constants $C_0,C_1>0$. 
\end{Lemma}
\begin{Lemma} Suppose $U$ is a plane domain, and has a compact subset $K$ and a non-empty open subset
$V$. Then there are positive constants $\delta_0$ and $C$ depending on $U$, $V$ and $K$, 
such that when $\delta<\delta_0$, the probability that a simple random walk on $\delta\Z^2$
 started from some $v\in\delta\Z^2\cap K$ will hit $V$ before exiting $U$ is greater than $C$.
\end{Lemma}
The following lemma about simple random walks on $D^\delta$ is an easy consequence of the above 
two lemmas and the Markov property of random walks. 
\begin{Lemma} For every $d>0$, there are $\delta_0, C>0$ depending on $d$ such that 
if $\delta<\delta_0$ and $v\in\delta\Z^2\cap D$ is such that $dist^\#(v,B_1)>d$, 
then the probability
that a simple random walk on $D^\delta$ started from $v$ hits $B_2$ before $B_1$ is at least $C$. 
\end{Lemma}
\begin{Lemma} For every $q\in(0,p)$ and $\eps>0$, there are $d,\delta_0>0$ depending on $q$ and 
$\eps$ such that for $\delta<\delta_0$, the probability that $dist^\#(\beta^{\delta}[q,p),B_1)\ge d$
is at least $1-\eps$. \end{Lemma}
{\bf Proof}. For $k=1,2,3$, let $J_k=W^{-1}(\CC_{q/k})$. Then $J_1$, $J_2$, $J_3$ are disjoint Jordan 
curves in $D$ that separate $B_2$ from $B_1$. And $J_2$ lies in the domain, denoted by 
$\Lambda$, bounded by $J_1$ and $J_3$. Moreover, the modulus of the domain bounded by $J_k$ 
and $B_2$ is $p-q/k$. Let $\tau^\delta$ be the first $n$ such that the edge $[y^{\delta}_{n-1}, y^{\delta}_{n}]$ intersects $J_2$. Then $\tau^\delta$ is a stopping time.
If $\delta$ is smaller than the distance between $J_1\cup J_3$ and $J_2$, then $y^{\delta}_{\tau^\delta}\in \Lambda$
and $y^{\delta}[-1,\tau^\delta]$ does not intersect $J_1$. Thus 
$M(D\sem y^{\delta}(-1,\tau^\delta])\ge p-q$, and so $T(\tau^\delta)\le q$.
So it suffices to prove that when $\delta$ and $d$ are small enough, the probability 
that $y^{\delta}$ will get within spherical distance $d$ from $B_1$ after time $\tau^\delta$ 
is less than $\eps$. Let $\RW^\delta$ denote a simple random walk on $D^\delta$ stopped on hitting $\partial D$, 
and $\CRW^\delta$ denote that $\RW^\delta$ conditioned to hit $B_2$ before $B_1$. Let $\RW^\delta_v$ and 
$\CRW^\delta_v$ denote that $\RW^\delta$ and $\CRW^\delta$, respectively, 
started from $v$. Since $y^{\delta}$ is obtained by erasing loops of CRW$_\delta^\delta$,
it suffices to show that the probability that CRW$_\delta^\delta$ will 
get within spherical distance $d$ from $B_1$ after it hits $\Lambda$, tends to zero as $d,\delta\to 0$. 
Since CRW$^\delta$ is a Markov chain, 
it suffices to prove that the probability that CRW$^\delta_v$ will
get within spherical distance $d$ from $B_1$ tends to zero as $d,\delta\to 0$, 
uniformly in $v\in\delta\Z^2\cap \Lambda$. By Lemma 3.9, there is $a>0$ such that for $\delta$ small enough, the probability that RW$^\delta_v$ hits $B_2$ before $B_1$ is greater than $a$,
for all $v\in\delta\Z^2\cap \Lambda$. 
By Markov property, for every $v\in\delta\Z^2\cap \Lambda$, the probability that 
CRW$^\delta_v$ will get within spherical distance $d$ from $B_1$ is 
less than $$\frac 1a \cdot\sup\{\PP[\RW^\delta_w\mbox{ hits }B_2\mbox{ before }B_1]:w\in V(D^\delta)
\cap D\mbox{ and }dist^\#(w,B_1)<d\},$$ which tends to $0$ as $d,\delta\to 0$ by Lemma 3.7. 
So the proof is finished. $\Box$

\begin{Lemma} For every $q\in(0,p)$ and $\eps>0$, there are $M,\delta_0>0$ depending on $q$ and 
$\eps$ such that for $\delta<\delta_0$, the probability that $\beta^{\delta}[q,p)\subset\B(0;M)$
is at least $1-\eps$. 
\end{Lemma}
{\bf Proof.} We use the notations of the last lemma. It suffices to prove that
the probability that $\RW^\delta_v\not\subset\B(0;M)$ tends to zero as $\delta\to 0$ and
 $M\to\infty$, uniformly in $v\in\delta\Z^2\cap \Lambda$. Let $K=\C\sem D$, then $K$ is 
unbounded, and the distance between $v\in\Lambda$ and $K$ is uniformly bounded from below by some 
$d>0$. Let $r>0$ be such that $\Lambda\subset\B(0;r)$. For $M>r$,
let $R=M-r$, then for $v\in\delta\Z^2\cap \Lambda$, $\RW^\delta_v$ should exit $\B(v;R)$ before $\B(0;M)$. By Lemma 3.7, the probability that $\RW^\delta_v\not\subset\B(0;M)$ is less than $C_0((\delta+d)/(M-r))^{C_1}$, which tends to $0$ as $\delta\to0$ and $M\to\infty$, uniformly in 
$v\in\delta\Z^2\cap \Lambda$. $\Box$

\begin{Lemma} For every $\eps>0$, there are $q\in(0,p)$ and $\delta_0>0$ depending on $\eps$ such that
when $\delta<\delta_0$, with probability greater than $1-\eps$, 
the diameter of $\beta^\delta[q,p)$ is less than $\eps$. \end{Lemma}
{\bf Proof.} The idea is as follows.
Note that as $q\to p$, the modulus of 
$D\sem \beta^\delta(0,q]$ tends to zero. So for any fixed $a\in (0,p)$, the spherical distance between 
$\beta^\delta[a,q]$ and $B_2$ tends to zero as $q\to p$. By Lemma 3.11, if $M$ is big and
$\delta$ is small, the fact that $\beta^\delta [a,q]$ does not
lie in $\B(0;M)$ is an event of small probability. Thus on the complement of this event, 
the Euclidean distance between $\beta^\delta[a,q]$ and $B_2$ tends
to zero, which means that $\beta^\delta$ gets to some point near $B_2$ in the Euclidean metric before time $q$. 
By Lemma 3.7, RW$^\delta_v$ does not go far before hitting $\partial D$ if $v$ is near $B_2$.
The same is true for CRW$^\delta_v$ because by Lemma 3.9, RW$^\delta_v$ hits $B_2$ 
before $B_1$ with a probability bigger than some positive constant when $v$ is near $B_2$.
Since $y^\delta$ is the loop-erasure of 
CRW$^\delta$, $y^\delta$ does not go far after it gets near $B_2$, nor does 
$\beta^\delta$. So the diameter of $\beta^\delta[q,p)$ is small. $\Box$

\begin{Definition} Let $z\in\C$, $r,\eps>0$. A $(z,r,\eps)$-quasi-loop in a path $\omega$ is a pair $a,b\in
\omega$ such that $a,b\in\mbox{\bf B}(z;r)$, $|a-b|\le\eps$, and the subarc of $\omega$ with endpoints $a$ and $b$
is not contained in $\mbox{\bf B}(z;2r)$. Let ${\cal L}^\delta(z,r,\eps)$ denote the 
event that $\beta^{\delta}[0,p)$ has a $(z,r,\eps)$-quasi-loop. \end{Definition}

\begin{Lemma} If $\lin{\mbox{\bf B}(z;2r)}\cap B_1=\emptyset$, then $\lim_{\eps\to 0}\mbox{\bf P}
[{\cal L}^\delta(z,r,\eps)]=0$, uniformly in $\delta$. \end{Lemma}
{\bf Proof}. This lemma is very similar to Lemma 3.4 in \cite{14}. There are two points of 
difference between them. First, here we are dealing with the loop-erased {\it conditional} random walk.
With Lemma 3.9, the hypothesis $\lin{\mbox{\bf B}(z;2r)}\cap B_1=\emptyset$ guarantees that
for some $v$ near $\partial\mbox{\bf B}(z;2r)$, 
 the probability that $\RW^\delta_v$ hits $B_2$ before
$B_1$ is bounded away from zero uniformly. Second, our LERW is stopped when it hits $B_2$, while in Lemma 
3.4 in \cite{14}, the LERW is stopped when it hits some single point. It turns out that the current setting is easier to deal with. See \cite{14} for more details. $\Box$

\begin{Proposition} For every $q\in(0,p)$ and $\eps>0$, there are $\delta_0,a_0>0$ 
depending on $q$ and
$\eps$ such that for $\delta<\delta_0$, with probability at least $1-\eps$, $\beta^\delta$ satisfies the following property. If $q\le t_1<t_2<p$, and 
$|\beta^\delta(t_1)-\beta^\delta(t_2)|<a_0$, then the diameter of 
$\beta^\delta[t_1,t_2]$ is less than $\eps$. \end{Proposition}
{\bf Proof}. For $d,M>0$, let $\Lambda_{d,M}$ denote the set of $z\in\B(0;M)$ such that $dist^\#(z,B_1)\ge d$, and
${\cal A}^\delta_{d,M}$ denote the event that $\beta^\delta[q,p)\subset\Lambda_{d,M}$. 
By Lemma 3.10 and 3.11, there are $d_0,M_0,\delta_0>0$ such that for $\delta<\delta_0$, 
$\PP[{\cal A}^\delta_{d_0,M_0}]>1-\eps/2$. 
Note that the Euclidean distance between $\Lambda_{d_0,M_0}$ and $B_1$ is greater than $d_0/2$. 
Choose $0<r<\min\{\eps/4,d_0/4\}$. 
There are finitely many points $z_1,\dots,z_n\in\Lambda_{d_0,M_0}$ such that $\Lambda_{d_0,M_0}\subset\cup_1^n\B(z_j;r/2)$.
For $a>0$, $1\le j\le n$, let ${\cal B}^\delta_{j,a}$ denote the event that $\beta^\delta[0,p)$ does not have a 
$(z_j,r,a)$-quasi-loop. Since $r<d_0/4$, we have $\lin{\mbox{\bf B}(z_j;2r)}\cap B_1=\emptyset$. 
By Lemma 3.13, there is $a_0\in(0,r/2)$ such that $\PP[{\cal B}^\delta_{j,a_0}]\ge 1-\eps/(2n)$ for $1\le j\le n$.
Let ${\cal C}^\delta=\cap_1^n{\cal B}^\delta_{j,a_0}\cap{\cal A}^\delta_{d_0,M_0}$. 
Then $\PP[{\cal C}^\delta]
> 1-\eps$ if $\delta<\delta_0$. And on the event ${\cal C}^\delta$, if there are $t_1<t_2\in[q,p)$ satisfying
$|\beta^\delta(t_1)-\beta^\delta(t_2)|<a_0$, then $\beta^\delta(t_1)$ lies in some ball $\mbox{\bf B}(z_j;r/2)$, so 
$\beta^\delta(t_2)\in\mbox{\bf B}(z_j;r)$ as $a_0<r/2$. Since $\beta^\delta$ does not have a $(z_j,r,a_0)$-quasi-loop, 
$\beta^\delta[t_1,t_2]\subset\mbox{\bf B}(z_j;2r)$. This then implies that the diameter of 
$\beta^\delta[t_1,t_2]$ is not bigger than $4r$, which is less than $\eps$. $\Box$
\vskip 2mm

Before the proof of Theorem 1.2, we need the notation of convergence of plane domain 
sequences. We say that a sequence of plane domains
$\{\Omega_n\}$ converges to a plane domain $\Omega$, or $\Omega_n\to \Omega$, if\\
(i) every compact subset of $\Omega$ lies in $\Omega_n$, for $n$ large enough;\\
(ii) for every $z\in\partial \Omega$ there exists $z_n\in\partial \Omega_n$ for each $n$ such that $z_n\to z$.\\
Note that a sequence of domains may have more than one limits. 
The following lemma is similar to Theorem 1.8, the Carath\'eodory kernel theorem, in \cite{11}. 
\begin{Lemma} Suppose $\Omega_n\to \Omega$, $f_n$ maps $\Omega_n$ conformally onto $G_n$, and $f_n$ converges to some function
$f$ on $\Omega$ uniformly on each compact subset of $\Omega$. Then either $f$ is constant on $\Omega$, or $f$ maps $\Omega$ conformally
onto some domain $G$. And in the latter case, $G_n\to G$ and $f_n^{-1}$ converges to $f^{-1}$ uniformly on
each compact subset of $G$. \end{Lemma}

\noindent
{\bf Proof of Theorem 1.2.} Suppose $(\chi^0_t,0\le t<p)$ is the driving function of 
$(W(K^0_t),0\le t<p)$. By Theorem 3.1, we may assume that all $\chi^\delta$ and $\chi^0$ are 
in the same probability space, so that for every $q\in(0,p)$ and $\eps>0$ there is an $\delta_0>0$ 
depending on $q$ and $\eps$ such that for $\delta<\delta_0$, 
$$\PP[\sup\{|\chi^\delta_t-\chi^0_t|:t\in [0,q]\}>\eps]<\eps.$$
Since $\beta^\delta$ and $\beta^0$ are determined by $\chi^\delta$ and $\chi^0$, respectively,
 all $\beta^\delta$ and $\beta^0$ are also in the same probability space. 
For the first part of this theorem, it suffices to prove that for every $q\in(0,p)$ 
and $\eps>0$ there is $\delta_0=\delta_0(q,\eps)>0$ such that for $\delta<\delta_0$,
\begin{equation}\PP[\sup\{|\beta^\delta(t)-\beta^0(t)|:t\in [q,p)\}>\eps]<\eps.\end{equation}

Now choose any sequence $\delta_n\to 0$. Then it contains a subsequence $\delta_{n_k}$ such that 
for each $q\in (0,p)$, $\chi^{\delta_{n_k}}$ converges to $\chi^0$ uniformly on $[0,q]$ almost surely. Here we use the fact that a sequence converging in probability contains an a.s.\ 
converging subsequence. For simplicity, we write $\delta_n$ instead of $\delta_{n_k}$. 
Let $\vphi^{\delta_n}_t$ ($\vphi^0_t$, resp.), $0\le t<p$, be 
the standard annulus LE maps of modulus $p$ driven by $\chi^{\delta_n}_t$
($\chi^0_t$, resp.), $0\le t<p$. Let $\Omega^{\delta_n}_t:=\A_p\sem W(\beta^{\delta_n}(0,t])$, 
and $\Omega^0_t:=\A_p\sem W(\beta^0(0,t])$. Fix $q\in(0,p)$. Suppose $K$ is a compact subset of $\Omega^0_q$. 
Then for every $z\in K$, $\vphi^0_t(z)$ does not blow up on $[0,q]$. Since the driving function 
$\chi^{\delta_n}$ converges to $\chi^0$ uniformly on $[0,q]$, so if $n$ is big enough, then for every $z\in K$,
$\vphi^{\delta_n}_t(z)$ does not blow up on $[0,q]$, which means that $K\subset \Omega^{\delta_n}_q$.
Moreover, $\vphi^{\delta_n}_q$ converges to $\vphi^0_q$ uniformly on $K$. It follows that $\Omega^{\delta_n}_q\cap
\Omega^0_q\to \Omega^0_q$ as $n\to\infty$. By Lemma 3.14, 
$(\vphi^{\delta_n}_q)^{-1}$ converges to $(\vphi^0_q)^{-1}$ uniformly on 
each compact subset of $\A_{p-q}$, and so 
$\Omega^{\delta_n}_q=(\vphi^{\delta_n}_q)^{-1}(\A_{p-q})\to
(\vphi^0_q)^{-1}(\A_{p-q})=\Omega^0_q$. Now we denote $D^{\delta_n}_t:=D\sem \beta^{\delta_n}(0,t]=W^{-1}(\Omega
^{\delta_n}_t)$, and $D^0_t:=D\sem \beta^0(0,t]=W^{-1}(\Omega^0_t)$. 
Then we have $D^{\delta_n}_q\to D^0_q$ for every $q\in (0,p)$. 

Fix $\eps>0$ and $q_1<q_2\in(0,p)$. Let $q_0=q_1/2$ and $q_3=(q_2+p)/2$. By Proposition 3.5, 
there are 
$n_1\in\N$ and $a\in(0,\eps/2)$ such that for $n\ge n_1$, with probability at least $1-\eps/3$, 
$\beta^{\delta_n}$ satisfies: if $q_0\le t_1<t_2<p$, and $|\beta^{\delta_n}(t_1)-\beta^{\delta_n}(t_2)|<a$, then the 
diameter of $\beta^{\delta_n}[t_1,t_2]$ is less than $\eps/3$. Let ${\cal A}_n$ denote the corresponding event.
Since $\beta^0$ is continuous, there is $b>0$ such that with probability $1-\eps/3$, we have $|\beta^0(t_1)-\beta^0(t_2)|< a/2$ if $t_1,t_2\in[q_0,q_3]$ and $|t_1-t_2|\le b$. Let ${\cal B}$ denote the 
corresponding event. We may choose $q_0<t_0<t_1=q_1<\cdots<t_{m-1}=q_2<t_m<q_3$
such that $t_{j}-t_{j-1}<b$ for $1\le j\le m$. Since $\beta^0(t_{j})\not\in \beta^0(0,t_{j-1}]$ for $1\le j\le m$,
there is $r\in(0,a/4)$ such that with probability at least $1-\eps/3$, $\lin{\mbox{\bf B}(\beta^0(t_{j});r)}
\subset D^0_{t_{j-1}}$ for all $0\le j\le m$. 
We now use the convergence of $D^{\delta_n}_t$ to $D^0_t$ for $t=t_0,\dots,t_m$. There exists $n_2\in\N$
such that for $n\ge n_2$, with probability at least $1-\eps/3$, $\lin{\mbox{\bf B}(\beta^0(t_{j});r)}
\subset D^{\delta_n}_{t_{j-1}}$, and there is some $z^n_{j}\in\partial D^{\delta_n}_{t_j}\cap\mbox{\bf B}
(\beta^0(t_{j});r)$, for all $1\le j\le m$. Let ${\cal C}_n$ denote the corresponding event. Then on the event
${\cal C}_n$, $z^n_j\in\partial D^{\delta_n}_j\sem\partial D^{\delta_n}_{j-1}$, so $z^n_j=\beta^{\delta_n}(s^n_j)$ 
for some $s^n_j\in(t_{j-1},t_j]$. Let ${\cal D}_n={\cal A}_n\cap{\cal B}\cap{\cal C}_n$. 
Then $\mbox{\bf P}[{\cal D}_n]\ge 1-\eps$, for $n\ge n_1+n_2$. And on the event ${\cal D}_n$, 
$$|z^n_j-z^n_{j+1}|\le 2r+|\beta^0(t_j)-\beta^0(t_{j+1})|\le 2r+a/2<a\mbox{, }\mbox{ }
\forall 1\le j\le m-1,$$
as $|t_j-t_{j+1}|\le b$.
Thus the diameter of $\beta^{\delta_n}[s^n_j,s^n_{j+1}]$ is less than $\eps/3$. 
It follows that for any $t\in [s^n_j,s^n_{j+1}]\subset[t_{j-1},t_{j+1}]$,
$$|\beta^0(t)-\beta^{\delta_n}(t)|\le|\beta^0(t)-\beta^0(t_j)|+|\beta^0(t_j)-z^n_j|
+|z^n_j-\beta^{\delta_n}(t)|\le a/2+r+\eps/3<\eps.$$
Since $[q_1,q_2]=[t_1,t_{m-1}]\subset\cup_{j=1}^{m-1}[s^n_j,s^n_{j+1}]$, we have now proved 
that for $n$ big enough,
with probability at least $1-\eps$, $|\beta^{\delta_n}(t)-\beta^0(t)|<\eps$ for all $t\in[q_1,q_2]$. By Lemma 3.12, 
for any $\eps>0$, there is $q(\eps)\in(0,p)$ such that if $n$ is big enough, with probability at least $1-\eps$,
the diameter of $\beta^{\delta_n}[q(\eps),p)$ is less than $\eps$. 
For any $S\in[q(\eps),p)$, by the uniform convergence
of $\beta^{\delta_n}$ to $\beta^0$ on the interval $[q(\eps),S]$, it follows that with
probability at least $1-\eps$, the
diameter of $\beta^0[q(\eps),S)$ is no more than $\eps$, nor is the diameter of $\beta^0[q(\eps),p)$. Now for fixed
$q\in(0,p)$ and $\eps>0$, choose $q_1\in(q,p)\cap(q(\eps/3),p)$. Then with probability at least $1-\eps/3$,
the diameter of $\beta^0[q_1,p)$ is less than $\eps/3$. And if $n$ is big enough, then with probability at least 
$1-\eps/3$, the diameter of $\beta^{\delta_n}[q_1,p)$ is less than $\eps/3$. Moreover, if $n$ is big enough, we may
require that with probability at least $1-\eps/3$,
$|\beta^{\delta_n}(t)-\beta^0(t)|\le\eps/3$ for all $t\in[q,q_1]$.
Thus $|\beta^{\delta_n}(t)-\beta^0(t)|\le\eps$ for all $t\in[q,p)$ with probability 
at least $1-\eps$, if $n$ is big
enough. Since $\{\delta_n\}$ is chosen arbitrarily, we proved formula (3.10). 

Now suppose that the impression of $0_+$ is the a single point, which must be $0$. From \cite{11}, we see that
$W^{-1}(z)\to 0$ as 
$z\in\A_p$ and $z\to 1$. From above, it suffices to prove that for any $\eps>0$, we can choose
$q\in(0,p)$ and $\delta_0>0$ such that for $\delta<\delta_0$, with probability at least $1-\eps$, the diameters of $\beta^\delta(0,q]$ and $\beta^0(0,q]$ are less than $\eps$. Since $W^{-1}$ is continuous at $1$, we need only to prove the
same is true for the diameters of $W(\beta^\delta(0,q])$ and $W(\beta^0(0,q])$. Note that they 
are the standard annulus LE 
hulls of modulus $p$ at time $q$, driven by $\chi^\delta_t$ and $\chi^0_t$, respectively. By 
Theorem 3.1, if $\delta$ and $q$ are small, then the diameters of $\chi^\delta[0,q]$ and $\chi^0[0,q]$ are uniformly
small with probability near $1$, so are the diameters of $W(\beta^\delta(0,q])$ and $W(\beta^0(0,q])$. $\Box$

\begin{Corollary} Almost surely $\lim_{t\to p}\beta^0(t)$ exists on $B_2$ . And the law 
is the same as the hitting point of a Brownian excursion in $D$
started from $0_+$ conditioned to hit $B_2$.
\end{Corollary}
A Brownian excursion in $D$ started from $0_+$ conditioned to hit $B_2$ is a random
closed subset of $D$ whose law is the weak limit as $\eps\to 0$ of the laws of Brownian 
motions in $D$ started from $\eps>0$ stopped on hitting $\pa D$ and conditioned to hit $B_2$.
\vskip 2mm
\noindent{\bf Proof of Corollary 1.1.} Now we consider the Riemann surface $R_p=(\R/(2\pi\Z))\times(0,p)$. Let $X_0=(\R/(2\pi\Z))\times\{0\}$ and $X_p=(\R/(2\pi\Z))\times\{p\}$ be the two boundary
components of $R_p$. Then $(x,y)\mapsto e^{-y+ix}$ is a conformal map from $R_p$ onto $\A_p$, and it maps $X_0$ and $X_p$
onto $\CC_0$ and $\CC_p$, respectively. So it suffices to prove this corollary with $\A_p$,
$\CC_p$ and $\CC_0$ replaced by $R_p$, $X_p$ and $X_0$, respectively.

For $n\in\N$, let $G_n$ be a graph that approximates $R_p$. The vertex set
$V(G_n)$ is 
$$\{(2k\pi/n,2m\pi/n):1\le k\le n, 0\le m\le\lfloor pn/(2\pi)\rfloor\}\cup\{(2k\pi/n,p):1\le k\le n\},$$
where $\lfloor x\rfloor$ is the maximal integer that is not bigger than $x$. And two vertices are connected by
an edge iff the distance between them is not bigger than $2\pi/n$. If $n>2\pi/p$, then for every vertex $v$ on 
$X_0$ or $X_p$, there is a unique $u\in V(G_n)\cap R_p$ that is adjacent to $v$. We write $u=N(v)$. For $v\in 
V(G_n)\cap X_0$, let $\RW$ be a simple random walk on $G_n$ started from $N(v)$ and stopped on hitting
$X_0\cup X_p$. Let $\CRW$ be that $\RW$ conditioned to hit $X_p$ before $X_0$.
Take the loop-erasure of $\CRW$, and then add the vertex $v$ at the beginning of the loop-erasure.
Then we get a simple lattice path from $v$ to $X_p$. We call this lattice path the LERW from $v$
to $X_p$. Similarly, for each $v\in X_p$, we may define the LERW from $v$ to $X_0$. Suppose 
$v\in V(D^\delta)\cap X_0$ and $u\in V(D^\delta)\cap X_p$. Let $P_{v,u}$ be the LERW from $v$ to $X_p$, conditioned 
to hit $u$, and $P_{u,v}$ be the LERW from $u$ to $X_0$, conditioned to hit $v$. By 
Lemma 7.2.1 in \cite{2}, the reversal 
of $P_{v,u}$ has the same law as $P_{u,v}$. Now we define the LERW from $X_0$ to $X_p$ to be the LERW
from a uniformly distributed random vertex on $X_0$ to $X_p$. Similarly, we may define the LERW from $X_p$ to $X_0$. 
It is clear that the hitting point at $X_p$ of the LERW from $X_0$ to $X_p$ is uniformly distributed. So the 
reversal of the LERW from $X_0$ to $X_p$ has the same law as the LERW from $X_p$ to $X_0$. 
Using the method in the proof of Theorem 1.2, we can show that 
the law of LERW from $X_0$ to $X_p$ converges to that of annulus SLE$_2$ in $R_p$
from a uniform random point on $X_0$ towards $X_p$. The same is true if we exchange 
$X_0$ with $X_p$. This ends the proof. $\Box$

\section{Disc SLE}
In this section, we will define another version of SLE: disc SLE, which describes a random process of growing 
compact subsets of a simply connected domain. Suppose $\Omega$ is a simply connected domain and $x\in \Omega$. Recall that a hull, say $F$, in $\Omega$ w.r.t.\ $x$, is a contractible compact subset of $\Omega$ that properly contains $x$. Then $\Omega\sem F$ is a doubly
connected domain with boundary components $\partial \Omega$ and $\partial F$. 
We say that $(F_t,a<t<b)$ is 
a Loewner chain in $\Omega$ w.r.t.\ $x$, if (i) each $F_t$ is a hull in $\Omega$ w.r.t.\ $x$; (ii) 
$F_s\subsetneqq F_t$ when $a<s<t<b$; and (iii) for any fixed $t_0\in(a,b)$, $(F_t\sem F_{t_0}, t_0\le t<b)$ is a Loewner chain in $\Omega\sem F_{t_0}$ on $\partial F_{t_0}$. 

\begin{Proposition} Suppose $\chi:(-\infty,0)\to\CC_0$ is continuous. Then there is a 
Loewner chain $(F_t,-\infty<t<0)$, in $\D$ w.r.t.\ $0$, and a family of maps $g_t$, 
$-\infty<t<0$, such that each $g_t$ maps $\D\sem F_t$ conformally onto
$\A_{|t|}$ with $g_t(\CC_0)=\CC_{|t|}$, and
\begin{equation}\left\{\begin{array}{lll}\partial_t g_t(z)=g_t(z)\SA_{|t|}(g_t(z)/\chi_t), &&-\infty<t<0;\\
   \lim_{t\to-\infty}e^t/g_t(z)=z, &&\forall z\in\D\sem\{0\}.\end{array}\right.\end{equation}
Moreover, such $F_t$ and $g_t$ are uniquely determined by $\chi_t$. 
We call $F_t$ and $g_t$, $-\infty<t<0$, the standard disc LE
hulls and maps, respectively, driven by $\chi_t$, $-\infty<t<0$.
\end{Proposition}

\noindent {\bf Proof.} For fixed $r\in(-\infty,0)$, let $\vphi^r_t$, $r\le t<0$, 
be the solution of
\begin{equation}\partial_t\vphi^r_t(z)=\vphi^r_t(z)\SA_{|t|}(\vphi^r_t(z)/\chi_t)\mbox{, }\mbox{ }\vphi^r_r(z)=z.\end{equation}
For $r\le t<0$, let $K^r_t$ be the set of $z\in\A_{|r|}$ such that $\vphi^r_s(z)$ blows up at some
time $s\in[r,t]$. Then $(K^r_t,r\le t<0)$ is a Loewner chain in 
$\A_{|r|}$ on $\CC_{0}$, 
and $\vphi^r_t$ maps $\A_{|r|}\sem K^r_t$ conformally onto $\A_{|t|}$ with $\vphi^r_t(\CC_{|r|})=\CC_{|t|}$.
By the uniqueness of the solution of ODE, if $t_1\le t_2\le t_3<0$, then
$\vphi^{t_2}_{t_3}\circ\vphi^{t_1}_{t_2}(z)=\vphi^{t_1}_{t_3}(z)$, for $z\in\A_{|t_1|}\sem K^{t_1}_{t_3}$. 
For $t<0$, define $R_t(z)=e^t/z$. Then $R_t$ maps $\A_{|t|}$ conformally onto itself, and exchanges the
two boundary components. Define $\ha{\vphi}^r_t=R_t\circ\vphi^r_t\circ R_r$, and $\ha{K}^r_t=R_r(K^r_t)$.
Then $\ha{K}^r_t$ is a hull in $\A_{|r|}$ on $\CC_{|r|}$, and $\ha{\vphi}^r_t$ maps $\A_{|r|}\sem\ha{K}^r_t$ 
conformally onto $\A_{|t|}$ with $\ha{\vphi}^r_t(\CC_0)=\CC_0$. We also have
$\ha{\vphi}^{t_2}_{t_3}\circ\ha{\vphi}^{t_1}_{t_2}(z)=\ha{\vphi}^{t_1}_{t_3}(z)$, for $z\in\A_{|t_1|}
\sem\ha{K}^{t_1}_{t_3}$, if $t_1\le t_2\le t_3<0$. And $\ha{\vphi}^r_t$ satisfies
$$\partial_t\ha{\vphi}^r_t(z)=\ha{\vphi}^r_t(z)\ha{\SA}_{|t|}(\ha{\vphi}^r_t(z)/\lin{\chi_t})\mbox{, }
\mbox{ }\ha{\vphi}^r_r(z)=z,$$
where $\ha{\SA}_p(z)=1-\SA_p(e^{-p}/z)$ for $p>0$. A simple computation gives:
$$|\ha{\SA}_p(z)|\le 8e^{-p}/|z|\mbox{, }\mbox{ }\mbox{if}\mbox{ }4e^{-p}\le |z|\le 1.$$ 
We then have 
\begin{equation} |\ha{\vphi}^r_t(z)-z|\le 8e^t\mbox{, }\mbox{ if }\mbox{ }r\le t<0\mbox{, }
\mbox{ and }\mbox{ }12e^t\le |z|\le 1.\end{equation} 

Now let $\ha{\psi}^r_t$ be the inverse of $\ha{\vphi}^r_t$. If $t_1\le t_2\le t_3$, then $\ha{\psi}^{t_1}_{t_2}\circ\ha{\psi}^{t_2}_{t_3}(z)=\ha{\psi}^{t_1}_{t_3}(z)$, for any $z\in\A_{|t_3|}$. 
For fixed $t\in(-\infty,0)$, $\{\ha{\psi}^r_t:r\in(-\infty,t]\}$ is a family of uniformly 
bounded conformal maps on $\A_{|t|}$, so is a normal family. This implies that we can find a sequence
$r_n\to-\infty$ such that for any $m\in\N$, $\{\ha{\psi}^{r_n}_{-m}\}$ converges to some $\ha{\psi}_{-m}$,
uniformly on each compact subset of $\A_m$. Let $\beta_n=\ha{\psi}^{r_n}_{-m}(\CC_{m/2})$. Then $\beta_n$ is 
a Jordan curve in $\A_{|r_n|}\sem\ha{K}^{r_n}_{-m}$ that separates the two boundary components. So $0$ is 
contained in the Jordan domain determined by $\beta_n$. Note that $\{\ha{\psi}^{r_n}_{-m}\}$ maps $\A_{m/2}$
onto the domain bounded by $\beta_n$ and $\CC_0$, whose modulus has to be $m/2$. So $\beta_n$ is not contained
in $\B(0;e^{-m/2})$. This implies that the diameter of $\beta_n$ is not less than $e^{-m/2}$. So 
$\ha{\psi}_{-m}$ can't be a constant. By Lemma 3.14, $\ha{\psi}_{-m}$ maps $\A_m$ conformally onto some domain 
$D_{-m}$, and $\ha{\psi}^{r_n}_{-m}(\A_m)\to D_{-m}$. Since $\ha{\psi}^{r_n}_{-m}(\A_m)=\A_{|r_n|}\sem\ha{K}
^{r_n}_{-m}\subset\D\sem\{0\}$, $D_{-m}\subset\D\sem\{0\}$. 
Since $M(\A_{|r_n|}\sem\ha{K}^{r_n}_{-m})=m$, there is some $a_m\in(0,1)$ such that
$\lin{\B(0;e^{r_n})}\cup\ha{K}^{r_n}_{-m}\subset\B(0;e^{-a_m})$ for all $r_n$.
So $\A_{a_m}$ contains no boundary points of $\A_{|r_n|}\sem\ha{K}^{r_n}_{-m}
=\ha{\psi}^{r_n}_{-m}(\A_m)$. Since these domains converge to $D_{-m}$ as $n\to\infty$,
so $\A_{a_m}$ contains no boundary points of $D_{-m}$, which means that either $\A_{a_m}
\subset D_{-m}$ or $\A_{a_m}\cap D_{-m}=\emptyset$. 
Now let $\gamma_n=\ha{\psi}^{r_n}_{-m}(\CC_{a_m/2})$. 
For the same reason as $\beta_n$, we have $\gamma_n\not\subset \B(0;e^{-a_m/2})$. 
So there is $z_n\in\CC_{a_m/2}$ such that
$|\ha{\psi}^{r_n}_{-m}(z_n)|\ge e^{-a_m/2}$. Let $z_0$ be any subsequential limit of $\{z_n\}$, then $z_0\in\CC_{a_m/2}\subset\A_m$ and $|\ha{\psi}_{-m}(z_0)|\ge e^{-a_m/2}$,
so $\ha{\psi}_{-m}(z_0)\in\A_{a_m}$. Thus $D_{-m}\cap \A_{a_m}\not=\emptyset$, and so $\A_{a_m}\subset D_{-m}$. Hence $D_{-m}$ has one boundary component $\CC_0$. Using similar arguments, we have $\ha{\psi}_t(\CC_0)=\CC_0$. 

If $r_n<-m_1<-m_2$, then $\ha{\psi}^{r_n}_{-m_1}\circ\ha{\psi}^{-m_1}_{-m_2}=\ha{\psi}^{r_n}_{-m_2}$,
which implies $\ha{\psi}_{-m_1}\circ\ha{\psi}^{-m_1}_{-m_2}=\ha{\psi}_{-m_2}$. 
For $t\in(-\infty,0)$, choose $m\in\N$ with $-m\le t$, define $\ha{\psi}_t=\ha{\psi}_{-m}
\circ\ha{\psi}^{-m}_t$ and $D_t=\ha{\psi}_t(\A_{|t|})$.
It is easy to check that the definition of $\ha{\psi}_t$ is independent of the choice of $m$, 
and the following properties hold. For all $t\in(-\infty,0)$,
$D_t$ is a doubly connected subdomain of $\D\sem\{0\}$ that has one boundary component $\CC_0$, and
$\ha\psi_t(\CC_0)=\CC_0$; $\ha{\psi}^{r_n}_t$ converges to $\ha{\psi}_t$, uniformly on each compact subset 
of $\A_{|t|}$. If $r<t<0$, then $\ha{\psi}_t=\ha{\psi}_r\circ\ha{\psi}^r_t$; 
$D_t\subsetneqq D_r$, and $D_r\sem D_t=\ha{\psi}_r(\ha{K}^r_t)$. 

Let $\ha{\vphi}_t$ on $D_t$ be the inverse of $\ha{\psi}_t$. By Lemma 3.14, $\ha{\vphi}^{r_n}_t$ converges to
$\ha{\vphi}_t$ as $n\to\infty$, uniformly on each compact subset of $D_t$. Thus 
from formula (4.3), we have
$|\ha{\vphi}_t(z)-z|\le 8e^t$, if
$12e^t\le |z|<1$. It follows that $\lim_{t\to-\infty}\ha{\vphi}_t(z)=z$, 
for any $z\in\D\sem\{0\}$. We also have
$\ha{\vphi}_t(z)=\ha{\vphi}^{-m}_t\circ\ha{\vphi}_{-m}(z)$, if $-m\le t<0$ and $z\in D_t$. 
Let $g_t=R_t\circ\ha{\vphi}_t$ on $D_t$. Then $g_t$ maps $D_t$ conformally onto $\A_{|t|}$,
takes $\CC_0$ to $\CC_{|t|}$, and
$$\lim_{t\to-\infty}e^t/g_t(z)=\lim_{t\to-\infty}\ha{\vphi}_t(z)=z\mbox{, }\mbox{ for any }z\in\D\sem\{0\}.$$
If $-m\le t$, then $g_t(z)=\vphi^{-m}_t\circ R_{-m}\circ\ha{\vphi}_{-m}(z)$, $\forall z\in D_t$. By formula (4.2), 
we have $$\partial_t g_t(z)=g_t(z)\SA_{|t|}(g_t(z)/\chi_t)\mbox{, }\mbox{ }-m\le t<0.$$
Since we may choose $m\in\N$ arbitrarily, formula (4.1) holds.

Let $F_t=\D\sem D_t$. Since $D_t$ is a doubly connected subdomain of $\D\sem\{0\}$ 
with a boundary component $\CC_0$, 
$F_t$ is a hull in $\D$ w.r.t.\ $0$. If $t_1<t_2<0$, then $F_{t_1}\subsetneqq F_{t_2}$, as $D_{t_1}\supsetneqq D_{t_2}$. 
Fix any $r\in(-\infty,0)$. For $t\in[r,0)$, $F_t\sem F_r=D_r\sem D_t=\ha{\psi}_r(\ha{K}^r_t)$. 
From Proposition 2.1 and the conformal invariance, $(\ha{\psi}_r(\ha{K}^r_t), r\le t<0)$ is a 
Loewner chain
in $D_r$ on $\partial F_r$. Thus $(F_t,-\infty<t<0)$ is a Loewner chain
in $\D$ w.r.t.\ $0$. 

Suppose $F_t^*$, $-\infty<t<0$, is a family of hulls in $\D$ on $0$, and $g_t^*$, $-\infty<t<0$, is a family
of maps such that for each $t$, $g_t^*$ maps $\D\sem F^*_t$ conformally onto $\A_{|t|}$ and  formula
(4.1) holds with $g_t$ replaced by $g_t^*$. By the uniqueness of the solution of ODE, we have $g_t^*=\vphi^r_t\circ
g^*_r$, if $r\le t<0$. So $R_t\circ g^*_t=\ha{\vphi}^r_t\circ R_r\circ g^*_r$. 
Now choose $r=r_n$ and let $n\to\infty$. Since 
$R_{r_n}\circ g_{r_n}^*\to\mbox{id}$ by formula (4.1) and $\ha{\vphi}^{r_n}_t\to\ha{\vphi}_t$, 
so $R_t\circ g^*_t=\ha{\vphi}_t$, from which follows that $g^*_t=R_t\circ\ha{\vphi}_t=g_t$ 
and $F^*_t= F_t$. $\Box$

\begin{Proposition} Suppose $(F_t,-\infty<t<0)$ is a Loewner chain in $\D$ w.r.t.\ $0$
such that $M(\D\sem F_t)=|t|$ for each $t$. Then there is a continuous $\chi:(-\infty,0)\to\CC_0$
 such that
$F_t$, $-\infty<t<0$, are the standard disc LE hulls driven by $\chi_t$, $-\infty<t<0$. \end{Proposition}
{\bf Proof.}  
For each $t<0$, choose $\vphi_t^*$ which maps $\D\sem F_t$ conformally onto $\A_{|t|}$ so that
$\vphi_t^*(1)=1$. Let $g_t^*=R_t\circ\vphi_t^*$, where $R_t(z)=e^t/z$. Then $g_t^*$ maps
$\D\sem F_t$ conformally onto $\A_{|t|}$ with $g_t^*(\CC_0)=\CC_{|t|}$ and $g_t^*(1)=e^t$. 
For any $r\le t<0$, let $K^*_{r,t}=g_r^*(F_t\sem F_r)$. Then for fixed $r<0$, $(K^*_{r,t},r\le t<0)$ is a Loewner chain in $\A_{|r|}$ on $\CC_0$. Now $g_t^*\circ(g_r^*)^{-1}$ maps
$\A_{|r|}\sem K^*_{r,t}$ conformally onto $\A_{|t|}$, and satisfies $g_t^*\circ(g_r^*)^{-1}(e^r)=e^t$.
From the proof of Proposition 2.1, there exists some continuous $\chi_{r,\cdot}^*:[r,0)\to\CC_0$ such that
for $r\le t<0$,
$$\partial_tg_t^*\circ(g_r^*)^{-1}(w)=g_t^*\circ(g_r^*)^{-1}(w)[\SA_{|t|}(g_t^*\circ(g_r^*)^{-1}(w)/
\chi^*_{r,t})-i\Imm \SA_{|t|}(e^t/\chi^*_{r,t})].$$
It then follows that
$$\partial_tg_t^*(z)=g_t^*(z)[\SA_{|t|}(g_t^*(z)/\chi^*_{r,t})-i\Imm \SA_{|t|}(e^t/\chi^*_{r,t})]
\mbox{, }\mbox{ }r\le t<0.$$
So $\chi^*_{r_1,t}=\chi^*_{r_2,t}$ if $r_1,r_2\le t$. We then have a continuous $\chi^*:(-\infty,0)\to\CC_0$,
such that $$\partial_tg_t^*(z)=g_t^*(z)[\SA_{|t|}(g_t^*(z)/\chi^*_t)-i\Imm \SA_{|t|}(e^t/\chi^*_t)]
\mbox{, }\mbox{ }-\infty\le t<0.$$
Consequently,
$$\partial_t\vphi_t^*(z)=\vphi_t^*(z)[\ha{\SA}_{|t|}(\vphi_t^*(z)/\lin{\chi^*_t})-i\Imm\ha{\SA}_{|t|}
(\chi^*_t)]\mbox{, }\mbox{ }-\infty\le t<0.$$
Since $|\ha{\SA}_{|t|}(z)|\le 8e^t$ when $4e^t\le|z|\le 1$, $|\Imm\ha{\SA}_{|t|}(\chi^*_t)|$ decays 
exponentially as $t\to-\infty$. Let $\theta(t)=\int_{-\infty}^t \Imm\ha{\SA}_{|s|}(\chi^*_s)ds$,
$\vphi_t(z)=e^{i\theta(t)}\vphi_t^*(z)$, and $\chi_t=e^{-i\theta(t)}\chi_t^*$. Then $\vphi_t$ maps 
$\D\sem F_t$ conformally onto $\A_{|t|}$ with $\vphi_t(\CC_0)=\CC_0$, and 
$$\partial_t\ln\vphi_t(z)=\partial_t\ln\vphi_t^*(z)+i\theta'(t)=\ha{\SA}_{|t|}(\vphi_t^*/\lin{\chi_t^*})
=\ha{\SA}_{|t|}(\vphi_t/\lin{\chi_t}).$$
Thus $\partial_t\vphi_t(z)=\vphi_t(z)\ha{\SA}_{|t|}(\vphi_t(z)/\lin{\chi_t})$. 
From the estimation of $\ha{\SA}_{|t|}$, we have 
$$|\vphi_t(z)-\vphi_r(z)|\le 8e^t\mbox{, }\mbox{ if }\mbox{ }12e^t\le|\vphi_r(z)|\le 1\mbox{, }\mbox{ and }
\mbox{ }r\le t<0.$$

Since $F_t$ contains $0$ and $M(\D\sem F_t)=|t|$, the diameter of $F_t$ tends to zero 
as $t\to-\infty$. Let $D_t=\D\sem F_t$. Then for
any sequence $t_n\to-\infty$, we have $D_{t_n}\to\D\sem\{0\}$. Since $\vphi_{t_n}$ is uniformly bounded,
there is a subsequence that converges to some function $\vphi$ on $\D\sem\{0\}$ uniformly on each compact
subset of $\D\sem\{0\}$. By checking the image of $\CC_1$ under $\vphi_{t_n}$ similarly as in the proof of Proposition 4.1, we see that $\vphi$ cannot be constant. So by Lemma 3.14, $\vphi$ maps $\D\sem\{0\}$ conformally onto 
some domain $D_0$ which is a subsequential limit of $\A_{|t_n|}=\vphi_{t_n}(D_{t_n})$. 
Since $t_n\to-\infty$, $D_0$ has to be $\D\sem\{0\}$ and so $\vphi(z)=\chi z$ for some $\chi\in\CC_0$. Now this
$\chi$ may depend on the subsequence of $\{t_n\}$. But we always have $\lim_{t\to-\infty}|\vphi_t(z)|=|z|$ for
any $z\in\D\sem\{0\}$. Now fix $z\in\D\sem\{0\}$, there is $s(z)<0$ such that when $r\le t<s(z)$, we have 
$12e^t\le|\vphi_r(z)|\le 1$. Therefore $|\vphi_t(z)-\vphi_r(z)|\le 8e^t$ for $r\le t<s(z)$. 
Thus $\lim_{t\to-\infty}\vphi_t(z)$ exists for every $z\in\D\sem\{0\}$. Since we have a sequence 
$t_n\to-\infty$ such that $\{\vphi_{t_n}\}$ 
converges pointwise to $z\mapsto \chi^* z$ on $\D\sem\{0\}$ for some $\chi^*\in\CC_0$, 
so $\lim_{t\to-\infty}\vphi_t(z)=\chi^*z$, for all $z\in\D\sem\{0\}$. 
Finally, let $g_t(z)=R_t\circ\vphi_t(z/\chi^*)$. Then $g_t$ maps $\D\sem F_t$ conformally onto $\A_{|t|}$, takes $\CC_0$ to $\CC_{|t|}$, and satisfies (4.1). $\Box$
\vskip 2mm
We still use $B(t)$ to denote a standard Brownian motion on $\R$ started from $0$. Let $\bf x$ be some uniform
random point on $\CC_0$, independent of $B(t)$. For $\kappa>0$ and 
$-\infty<t<0$, write $\chi^\kappa_t={\bf x}e^{iB(\kappa|t|)}$. 
The process $(\chi^\kappa)$ is determined by the following properties: for any fixed
$r<0$, $(\chi^\kappa_t/\chi^\kappa_r,r\le t<0)$ has the same law as 
$(e^{iB(\kappa(t-r))},r\le t<0)$ and is independent from $\chi^\kappa_r$. 
If $F_t$ and $g_t$, $-\infty<t<0$, are the standard disc LE hulls and maps, respectively,
driven by $\chi^\kappa_t$, $-\infty<t<0$, then we call them the standard disc SLE$_\kappa$ hulls and
maps, respectively. From the properties of $\chi^\kappa_t$, we see that for any fixed $r<0$, $g_r(F_{r+t}\sem F_r)$, 
$0\le t<|r|$, is an annulus SLE$_\kappa(\A_{|r|};\chi^\kappa_r\to\CC_{|r|})$. 
The existence of standard annulus SLE$_\kappa$ trace then implies the a.s.\ existence 
of standard disc SLE$_\kappa$ trace, which is a curve $\gamma:[-\infty,0)\to\D$ such that $\gamma(-\infty)=0$, 
and for each $t\in(-\infty,0)$, $F_t$ is the hull generated by $\gamma[-\infty,t]$, i.e., the complement of the 
unbounded component of $\C\sem\gamma[-\infty,t]$. If $\kappa\le 4$, the trace is a simple curve; otherwise, 
it is not simple.
Suppose $D$ is a simply connected domain and $a\in D$. Let $f$ map $\D$ conformally
onto $D$ so that $f(0)=a$ and $f'(0)>0$. Then we define $f(F_t)$ and $f(\gamma(t))$, $-\infty\le t<0$, to be the disc SLE$_\kappa(D;a\to\partial D)$ hulls and trace. 

The next theorem is about the equivalence of disc SLE$_6$ and full plane SLE$_6$. First, let's review
the definition of full plane SLE. It was proved in \cite{10} that for any continuous $\chi:(-\infty,+\infty)\to\CC_0$, 
there is a Loewner chain $(F_t,-\infty<t<+\infty)$, in $\C$ w.r.t.\ $0$, and a family of 
maps $g_t$, $-\infty<t<+\infty$, such that for each $t$, $g_t$ maps $\ha{\C}\sem F_t$
conformally onto $\D$ with $g_t(\infty)=0$, and
\begin{equation*}\left\{\begin{array}{lll}\partial_t g_t(z)=g_t(z)\frac{1+g_t(z)/\chi_t}{1-g_t(z)/\chi_t}, 
&&-\infty<t<+\infty;\\ \lim_{t\to-\infty}e^t/g_t(z)=z, &&\forall z\in\C\sem\{0\}.\end{array}\right.\end{equation*}
Such $F_t$ and $g_t$, $-\infty<t<+\infty$, are unique, and are called the full plane LE hulls 
and maps, respectively, driven by $\chi_t$, $-\infty<t<+\infty$. 
The diameter of $F_t$ tends to $0$ as $t\to-\infty$; and tends to $\infty$ as $t\to+\infty$. 

The driving process of full plane SLE$_\kappa$ is an extension of $\chi^\kappa_t$ to $\R$ defined as follows.
Choose another standard Brownian motion $B'(t)$ on $\R$ started from $0$, which is independent
of $B(t)$ and $\bf x$. For $t\ge 0$, let $\chi^\kappa_t={\bf x}e^{iB'(t)}$. Then for any fixed $r\in\R$, 
$\chi^\kappa_t/\chi^\kappa_r$, $r\le t<+\infty$, have the same distribution as $e^{iB(\kappa(t-r))}$, 
$r\le t<+\infty$. This implies that for full plane SLE$_\kappa$ hulls $F_t$, $t\in\R$, and
any fixed $r\in\R$, $(g_r(F_{r+t}\sem F_r))$ has the same law as radial 
SLE$_\kappa(\D;\chi^\kappa_r\to 0)$. 

Suppose $\Omega$ is a simply connected plane domain that contains $0$. Let $\tau$ be the first
$t$ such that full plane SLE$_\kappa$ hull $F_t\not\subset \Omega$. Then as $t\nearrow \tau$, $F_t$ approaches $\partial \Omega$, and 
$(F_t, -\infty<t<\tau)$ is a Loewner chain in $\Omega$ w.r.t.\ $0$. Let $u(t)=-M(\Omega\sem F_t)$,
for $-\infty<t<\tau$. Then $u$ is a continuous increasing function, and maps $(-\infty,\tau)$ onto $(-\infty,0)$.
Let $v$ be the inverse of $u$, and choose $f$ that maps $\D$ onto $\Omega$ with $f(0)=0$ 
and $f'(0)>0$. Then 
$f^{-1}(F_{v(s)})$, $-\infty<s<0$, are the standard disc LE hulls driven by some function. 
Using the same method in the proof of Theorem 1.1, we can prove that this driving function
has the same law as $(\chi^6_t)_{-\infty<t<0}$. So we have

\begin{Theorem} Suppose $\Omega$ is a simply connected domain that contains $0$. Let 
$(K_t,-\infty<t<+\infty)$ be full plane SLE$_6$ hulls, and $(L_s,-\infty<s<0)$ be the
disc SLE$_6(\Omega;0\to\partial \Omega)$. Let $\tau$ be the first $t$ that $K_t\not\subset
\Omega$. Then up to a time-change, $(K_t,-\infty<t<\tau)$ has the same law as $L_s$,
$-\infty<s<0$. \end{Theorem}

\begin{Corollary} The distribution of the hitting point of full plane SLE$_6$ trace at 
$\partial \Omega$ is the harmonic measure valued at $0$. 
\end{Corollary}
An immediate consequence of this corollary is that the plane SLE$_6$ hull stopped at 
the hitting time 
of $\partial \Omega$ has the same law as the hull generated by a
plane Brownian motion started from $0$ and stopped on exiting $\Omega$. See \cite{17} and \cite{7} for details. 
\vskip 2mm
Disc SLE$_2$ is also interesting. Suppose $\Omega$ is a simply connected domain that contains $0$. 
Let $\RW$ be a simple random walk on $\Omega^\delta$ started from $0$, and stopped on hitting $\partial \Omega$. 
Let LERW be the loop-erasure of $\RW$. Then LERW is a simple lattice path from
$0$ to $\partial \Omega$. Write LERW as $y=(y_0,\dots,y_{\up})$ with $y_0=0$ and 
$y_{\up}\in\partial \Omega$. We may extend $y$ to be defined on $[0,\up]$ so that it 
is linear on each $[j-1,j]$ for $1\le j\le\up$. 
Then it is clear that $(y(0,s],0\le s<\up)$ is a Loewner chain in $\Omega$ w.r.t.\ $0$. 
Let $T(s)=-M(\Omega\sem y(0,s])$, for $0<s<\up$. Then $T$ is a continuous increasing function, 
and maps 
$(0,\up)$ onto $(-\infty,0)$. Let $S$ be the inverse of $T$. Define $\beta^\delta(t)=y(S(t))$, for 
$-\infty<t<0$, and $\beta^\delta(-\infty)=0$. 
Let $\beta^0:[-\infty,0)\to \Omega$ be the trace of disc SLE$_2(\Omega;0\to\partial \Omega)$. 

\begin{Theorem} For any $\eps>0$, there is $\delta_0>0$ such that for $\delta<\delta_0$, we may couple
$\beta^\delta$ with $\beta^0$ so that 
$$\PP[\sup\{|\beta^\delta(t)-\beta^0(t)|:-\infty\le t<0\}\ge \eps]<\eps.$$
\end{Theorem}
{\bf Proof.} Note that $\Omega^\delta$ may not be connected, we replace it by its 
connected component
that contains $0$. Let $g_0$ be constant $1$ on $V(\Omega^\delta)$. For $0<j<\up_\delta$, 
let $g_j$ be the $g$ in Lemma 3.4 with $A=V(\Omega^\delta)\cap\partial \Omega$, 
$B=\{y_0,\dots,y_{j-1}\}$, and $x=y_j$. 
Similarly as Proposition 3.2 and 3.3, $g_j$'s are observables for the LERW here, and
they approximate the observables for disc SLE$_2$. We may 
follow the process in proving Theorem 1.2. $\Box$

\begin{Corollary} Suppose $\Omega$ is a simply connected plane domain, and $a\in \Omega$. Let 
$\beta(s)$, $-\infty<s<0$,
be the disc SLE$_2(\Omega;a\to\partial \Omega)$ trace. Let $\gamma(t)$, $0< t<\infty$, be the 
radial SLE$_2(\Omega;{\bf x}\to 0)$ trace, where $\bf x$ is a random point 
on $\partial \Omega$ with harmonic measure at $a$. Then the reversal of $\beta$
has the same law as $\gamma$, up to a time-change. 
\end{Corollary}
{\bf Proof.} This follows immediately from Theorem 4.2, the 
approximation of LERW to radial SLE$_2$ in \cite{6}, and the reversibility property
of LERW in \cite{2}. $\Box$  

\section{Convergence of the observables}
This is the last section of this paper. The goal is to prove Proposition 3.3.
The proof is sort of long. The main difficulty is that we need the approximation to be uniform 
in the domains. The tool we can use is Lemma 3.14. However, the limit of a domain sequence
in general does not have good boundary conditions, even if every domain in the sequence has.
Prime ends and crosscuts are used to describe the boundary correspondence under conformal maps. 
Some ideas of the proof come from \cite{6}.

We will often deal with a function defined on a subset of $\delta\Z^2$. Suppose
$f$ is such a function. For $v\in\delta\Z^2$ and $z\in\Z^2$, 
if $f(v)$ and $f(v+\delta z)$ are defined, then define
$$\nabla^\delta_zf(v)=(f(v+\delta z)-f(v))/\delta,$$
We say that $f$ is $\delta$-harmonic in $\Omega\subset\C$
if $f$ is defined on $\delta\Z^2\cap \Omega$ and all $v\in\delta \Z^2$ that are adjacent to 
vertices of $\delta\Z^2\cap \Omega$ so that for all $v\in\delta\Z^2\cap \Omega$,
$$f(v+\delta)+f(v-\delta)+f(v+i\delta)+f(v-i\delta)=4f(v).$$
The following lemma is well known.
\begin{Lemma} Suppose $\Omega$ is a plane domain that has a compact subset $K$. 
For $l\in\N$, let $z_1,\dots,z_l\in\Z^2$. 
Then there are positive constants $\delta_0$ and $C$ depending on $\Omega$, $K$, and $z_1,\dots,z_l$,
such that for $\delta<\delta_0$, if $f$ is non-negative and $\delta$-harmonic in $\Omega$, then
for all $v_1,v_2\in\delta\Z^2\cap K$, $$\nabla^\delta_{z_1}\cdots\nabla^\delta_{z_l}f(v_1)\le Cf(v_2).$$
This is also true for $l=0$, which means that $f(v_1)\le C f(v_2)$. \end{Lemma}

For $a,b\in\delta\Z$, denote
$$S_{a,b}^\delta:=\{(x,y):a\le x\le a+\delta,b\le y\le b+\delta\}.$$ 
Suppose $A$ is a subset of $\delta\Z^2$, let $S_A^\delta$ be the union of 
all $S_{a,b}^\delta$ whose four vertices
are in $A$. If $f$ is defined on $A$, we may define a continuous function $\CE^\delta f$ 
on $S_A^\delta$, as follows. For $(x,y)\in S_{a,b}^\delta\subset S_A^\delta$, define
\begin{eqnarray*}\CE^\delta f(x,y)=&(1-s)(1-t)f(a,b)+(1-s)tf(a,b+\delta)\\
&+s(1-t)f(a+\delta,b)+stf(a+\delta,b+\delta),
\end{eqnarray*}
where $s=(x-a)/\delta$ and $t=(y-b)/\delta$. Then $\CE^\delta f$ is well defined on 
$S_A^\delta$, and agrees with $f$ on $S_A^\delta\cap A$. Moreover, on $S_{a,b}^\delta$, 
$\CE^\delta f$ has a Lipschitz constant
not bigger than two times the maximum of 
$|\nabla^\delta_{(1,0)}f(a,b)|$, $|\nabla^\delta_{(0,1)}f(a,b)|$,
$|\nabla^\delta_{(1,0)}f(a,b+\delta)|$, $|\nabla^\delta_{(0,1)}f(a+\delta,b)|$.
And for any $u\in\Z^2$, 
$$\CE^\delta\nabla_u^\delta f(z)=(\CE^\delta f(z+\delta u)
-\CE^\delta f(z))/\delta,$$ when both sides are defined.
\vskip 2mm 
\noindent
{\bf Proof of Proposition 3.3.} Suppose the proposition is not true. Then we can find $\eps_0>0$,
a sequence of lattice paths $w_n\in L^{\delta_n}$ with $\delta_n\to 0$, and a sequence of
points $v_n\in V^{\delta_n}_2$, such that $|g_{w_n}(v_n)-u_{w_n}(v_n)|>\eps_0$ for all $n\in\N$.
For simplicity of notations, we write $g_n$ for $g_{w_n}$, $u_n$ for $u_{w_n}$, and $D_n$ for $D_{w_n}$. 
Let $p_n$ be the modulus of $D_n$.  The remaining of the proof is composed of four steps.

\subsection{The limits of domains and functions}
By comparison principle of extremal length, we have $p\ge p_n\ge M(U_2)>0$.
By passing to a subsequence, we may assume that $p_n\to p_0\in(0,p]$. Then $\A_{p_n}\to \A_{p_0}$. 
Let $Q_n$ map $D_n$ conformally onto $\A_{p_n}$ so that $Q_n(z)\to 1$ as $z\in D_n$ and $z\to P(w_n)$. 
Then $u_n=\Ree\SA_{p_n}\circ Q_n$.
Now $Q_n^{-1}$ maps $\A_{p_n}$ conformally
onto $D_n\subset D$. Thus $\{Q_n^{-1}\}$ is a normal family. By passing to a subsequence, we may 
assume
that $Q_n^{-1}$ converges to some function $J$ uniformly on each compact subset of $\A_{p_0}$. 
Using some argument similar to that in the proof of Theorem 1.2, we conclude that
$J$ maps $\A_{p_0}$ conformally onto some domain $D_0$, and $D_n\to D_0$.  
Let $Q_0=J^{-1}$ and $u_0=\Ree \SA_{p_0}\circ Q_0$.
Then $Q_n$ and $u_n$ converge to $Q_0$ and $u_0$, respectively, uniformly on each compact subset of $D_0$. Moreover, we have $U_2\cup\alpha_2\subset D_0\subset D$. Thus 
$B_2$ is one boundary component of $D_0$. Let $B_1^n$ and $B_1^0$ denote the boundary component of $D_n$ and $D_0$, respectively, other than $B_2$.

 
Let $\{K_m\}$ be a sequence of compact subsets of $D_0$ such that $D_0=\cup_m K_m$, and for each $m$, $K_m$ disconnects $B_1^0$ from $B_2$ and $K_m\subset\mbox{int}\,K_{m+1}$. 
Let $K^n_m=K_m\cap\delta_n\Z^2$.
Now fix $m$. If $n$ is big enough depending on $m$, we can have the following properties. First, $K_m\subset D_n$ and $K_m^n\subset V(D^{\delta_n})$, so $g_n$ is $\delta_n$-harmonic on $K_m$. Second, $K_m^n$ disconnects all lattice 
paths on $D^{\delta_n}$ from $B_2$ to $B_1^n$. Now let $\RW_v^n$ be a simple random walk on $D^{\delta_n}$ started from 
$v\in V(D^{\delta_n})$, and $\tau^n_m$ the hitting time of $\RW_v^n$ at $B_2\cup K^n_m$. By the properties of $g_n$, 
if $v$ is in $D$ and between $K_m$ and $B_2$, then $(g_n(\RW_v^n(j)),0\le j\le\tau_n^m)$ is a martingale,
so $g_n(v)=\EE[g_n(\RW^{v}_n(\tau^m_n)]$. Now suppose $g_n(v)>1$ for all $v\in K^m_n$.  
Choose $v_0\in V(D^{\delta_n})\cap D$ that is adjacent to some vertex of $F^{\delta_n}=V(D^{\delta_n})\cap B_2$. 
Then $g_n(v_0)=\EE[g_n(\RW^{v_0}_n(\tau^m_n))]\ge 1$. The equality holds iff there is no lattice path on $D^{\delta_n}$ 
from $v_0$ to $K_m^n$. By the definition of $D^{\delta_n}$, we know that the equality can not always hold.
It follows that $\sum_{u\in F^{\delta_n}}\Delta_{D^{\delta_n}}
g_n(u)>0$, which contradicts the definition of $g_n$. Thus there is $v\in K_m^n$ such that $g_n(v)\le 1$. 
Note that $g_n$ is non-negative. By Lemma 5.1, if $n$ is big enough depending on $m$, then $g_n$ on $K_m^n$ is 
uniformly bounded in $n$. Similarly for any $z_1,\dots,z_l\in\Z^2$, $\nabla_{z_1}^{\delta_n}\cdots\nabla_{z_l}^{\delta_n}g_n$ 
on $K_m^n$ is uniformly bounded in $n$, if $n$ is big enough depending on $m$, and $z_1,\dots,z_l\in\Z^2$.

We just proved that for a fixed $m$, if $n$ is big enough depending on $m$, then $g_n$ on $K_{m+1}^n$ is 
$\delta_n$-harmonic and uniformly bounded in $n$. We may also choose $n$ big such that every lattice square of $\delta_n\Z^2$ that intersects $K_{m}$ is contained in $K_{m+1}$, and so $\CE^{\delta_n}g_n$ on $K_{m}$ is well 
defined, and is uniformly bounded in $n$. Using the boundedness of $\nabla_u^{\delta_n} g_n$ on $K_{m+1}^n$
for $u\in\{1,i\}$, we conclude that $\{\CE^{\delta_n}g_n\}$ on $K_{m}$ is uniformly continuous.
By Arzela-Ascoli Theorem, there is a 
subsequence of $\{\CE^{\delta_n}g_n\}$, which converges uniformly on $K_m$. By passing to a subsequence, we may assume 
that $\CE^{\delta_n}g_n$ converges uniformly on each $K_m$. Let $g_0$ on $D_0$ be the limit function. Similarly, for 
any $z_1,\dots,z_l\in\Z^2$, there is a subsequence of $\{\CE^{\delta_n}\nabla_{z_1}^{\delta_n}\cdots\nabla_{z_l}
^{\delta_n}g_n\}$ which converges uniformly on each $K_m$. By passing to a subsequence again, we may assume
that for any $z_1,\dots,z_l\in\Z^2$, $\CE^{\delta_n}\nabla_{z_1}^{\delta_n}\cdots\nabla_{z_l}^{\delta_n}g_n$
converges to $g_0^{z_1,\dots,z_l}$ on $D_0$, uniformly on each $K_m$. It is easy to check that $$g_0^{z_1,\dots,z_l}=(a_1\partial_x+b_1\partial_y)\cdots(a_l\partial_x+b_l\partial_y)g_0,$$ 
if $z_j=(a_j,b_j)$, $1\le j\le l$. Since $g_n$ is $\delta_n$-harmonic on $K_m$ for $n$ big enough, we have $(\nabla^{\delta_n}_1\nabla^{\delta_n}_{-1}+\nabla^{\delta_n}_i\nabla^{\delta_n}_{-i})g_n\equiv 0$ on $K_m^n$. 
Thus $(\partial_x^2+\partial_y^2)g_0=0$, which means that $g_0$ is harmonic. 

Now suppose $x_n\in V(D^{\delta_n})\cap D\to B_2$ in the spherical metric. Since the spherical distance 
between $K_1$ and $B_2$ is positive, the probability that a simple random walk on $D^{\delta_n}$ started
 from $x_n$ hits $K_1$ before $B_2$ tends to zero by Lemma 3.7. If $n$ is big enough, $K_1$ is a subset of $D_n$ and 
disconnects $B_2$ from $B_1^n$. We have proved that $g_n$ is uniformly bounded on 
$\delta_n\Z^2\cap K_1$, if $n$ is big enough. And by definition $g_n\equiv 1$ 
on $V(D^{\delta_n})\cap B_2$. By Markov
property, we have $g_n(x_n)\to 1$. Since $g_0$ is the limit of $\CE^{\delta_n} g_n$, this implies that 
$g_0(z)\to 1$ as $z\in D_0$ and $z\to B_2$ in the spherical metric. 
Thus $g_0\circ J(z)\to 1$ as $z\in \A_{p_0}$ and $z\to\CC_{p_0}$.

Now let us consider the behavior of $u_n$ and $u_0$ near $B_2$. If $z\in D_n$ and $z\to B_2$ 
in the spherical metric, 
then $Q_n(z)\to \CC_{p_n}$, and so $u_n(z)=\Ree \SA_{p_n}\circ Q_n(z)\to 1$. 
Using a plane Brownian motion instead of a simple random walk in the above argument, 
we conclude that $u_n(z)\to 1$ as $z\in D_n$ and $z\to B_2$ in the 
spherical metric, uniformly in $n$. 


Suppose $\{v_n\}$, chosen at the beginning of this proof, has a subsequence that tends to $B_2$ in the spherical 
metric. By passing to a subsequence, we may assume that $v_n\to B_2$ in the spherical metric.  From the result of the 
last two paragraphs, we see that $g_n(v_n)\to 1$ and $u_n(v_n)\to 1$. This contradicts the hypothesis that 
$|g_n(v_n)-u_n(v_n)|\ge\eps_0$. Thus $\{v_n\}$ has a positive spherical distance from $B_2$. 
Since the domain bounded by $\alpha_1$ and $\alpha_2$ disconnects $U_2$ from $B_1^0$, and
$\{v_n\}\subset U_2$, so $\{v_n\}$ has a positive spherical distance from $B_1$ too.
Thus $\{v_n\}$ has a 
subsequence that converges to some $z_0\in D_0$. Again we may assume that $v_n\to z_0$. Then 
$u_0(z_0)=\lim u_n(v_n)$ and $g_0(z_0)=\lim g_n(v_n)$, and so $|u_0(z_0)-g_0(z_0)|\ge\eps_0$. 
We will get a contradiction by proving that $g_0\equiv u_0$ in $D_0$.

Note that $g_0$ is non-negative, since
each $g_n$ is non-negative. We can find a Jordan curve $\beta$ in $D_0$ which satisfies the following properties. 
It disconnects $B_2$ from $B_1^0$; it is the union of finite line segments which are parallel to either $x$ or $y$ 
axis; and it does not intersect $\cup_n\delta_n\Z^2$. By Remark 2 in Section 3 and the uniform convergence of $\nabla^{\delta_n}_1 g_n$ to $\partial_x g_0$, and $\nabla^{\delta_n}_i g_n$ to $\partial_y g_0$ on some 
neighborhood of $\beta$, we have $\int_{\beta}\partial_{\bf n} g_0ds=0$, where $\mbox{\bf n}$ 
is the unit norm vector on $\beta$ pointed towards $B_1$. 
Thus $g_0$ has a harmonic conjugate, and so does $g_0\circ J$. 
We will prove $g_0\circ J=\Ree \SA_{p_0}$, from which follows that $g_0=u_0$. 
We have proved that $g_0\circ J(z)\to 1$ as $\A_{p_0}\ni z\to\CC_{p_0}$. It suffices to 
show that $g_0\circ J(z)\to 0$ as $\A_{p_0}\sem U\ni z\to \CC_{0}$ for any neighborhood
$U$ of $1$. 

\vskip 4mm
\subsection{The existence of some sequences of crosscuts}

For a doubly connected domain $\Omega$ and one of its boundary component $X$, 
we say that $\gamma$ is a crosscut in $\Omega$ on $X$ if $\gamma$ is an open simple curve 
in $D$ whose two ends approach two points (need not be distinct) of $X$ in Euclidean distance. 
For such $\gamma$, $\Omega\sem\gamma$ has two connected components, one is a simply connected 
domain, and the other is a doubly connected domain. Let $U(\gamma)$ denote the simply
connected component of $D\sem\gamma$. Then $\partial U(\gamma)$ is the union of $\gamma$
and a subset of $X$. 

Now $Q_0$ maps $D_0$ conformally onto $\A_{p_0}$, and $Q_0(B_1^0)=\CC_0$. 
Similarly as Theorem 2.15 in \cite{11}, we can find a sequence of 
crosscuts $\{\gamma^k\}$ in $D_0$ on $B_1^0$ which satisfies\\
(i) for each $k$, $\lin{\gamma^{k+1}}\cap\lin{\gamma^k}=\emptyset$ and 
$U(\gamma^{k+1})\subset U(\gamma^k)$;\\
(ii) $Q_0(\gamma^k)$, $k\in\N$, are mutually disjoint crosscuts in $\A_{p_0}$ on $\CC_0$; and\\
(iii) $U(Q_0(\gamma^k))$, $k\in\N$, forms a neighborhood basis of $1$ 
in $\A_{p_0}$.\\
Note that $U(Q_0(\gamma^k))=Q_0(U(\gamma^k))$, so $U(Q_0(\gamma^{k+1}))\subset
U(Q_0(\gamma^k))$, for all $k\in\N$.
We will prove that there is some crosscut $\gamma^k_n$ in each $D_n$ on $B_1^n$ such that 
$\gamma^k_n$ and $Q_n(\gamma^k_n)$ converge to $\gamma^k$ and 
$Q_0(\gamma^k)$, respectively, in the sense that we will specify.


Now fix $k\in\N$ and $\eps>0$.  
Parameterize $\lin{\gamma^k}$ and $\lin{Q_0(\gamma^k)}$ as the image of the function $a:[0,1]\to
D\cup B_1^0$ and $b:[0,1]\to \A_{p_0}\cup \CC_0$, respectively, so that $b(t)=Q_0
(a(t))$, for $t\in (0,1)$. We may choose $s_1\in(0,1/2)$ such that the diameters of $a[0,s_1]$ and $a[1-s_1,1]$
are both less than $\eps/3$. There is $r_1\in(0,\eps)\cap(0,(1-e^{-p_0})/2)$
such that the curve
$b[s_1,1-s_1]$ and the balls $\lin{\B(b(0);r_1)}$ and $\lin{\B(b(1);r_1)}$ are mutually disjoint.
Suppose $\gamma^k$ is contained in $\mbox{\bf B}(0;M)$ for some $M>\eps$. There is $C_M>0$ 
such that the spherical 
distance between any $z_1,z_2\in \mbox{\bf B}(0;2M)$ is at least $C_M|z_1-z_2|$. 
So for every smooth curve $\gamma$ in $\B(0;2M)$, we have $L^\#(\gamma)\ge C_ML(\gamma)$, where
$L$ and $L^\#$ denote the Euclidean length and spherical length, respectively. 
Let $r_2=r_1\exp(-72\pi^2/(C_M^2\eps^2))$. Then we may choose $s_2\in(0,s_1)$ such that $b[0,s_2]\subset 
\mbox{\bf B}(b(0);r_2)$ and $b[1-s_2,1]\subset \mbox{\bf B}(b(1);r_2)$. 

For $j=0,1$, let $\Gamma_j$ be the set of crosscuts $\gamma$ in $\A_{p_0}$ on $\CC_0$ such
that $$\mbox{\bf B}(b(j);r_2)\cap \D\subset U(\gamma)\subset \mbox{\bf B}(b(j);r_1).$$ 
Then the extremal length of $\Gamma_j$ is less than 
$$2\pi/(\ln r_1-\ln r_2)=C_M^2\eps^2/(36\pi).$$ 
If $n$ is big enough, then
$\mbox{\bf B}(b(j);r_1)\cap \D\subset \A_{p_n}$, so all $\gamma\in \Gamma_j$ are in $\A_{p_n}$. Then the
extremal length of $Q_n^{-1}(\Gamma_j)$ is also less than $C_M^2\eps^2/(36\pi)$. 
Since the spherical area of $Q_n^{-1}(\A_{p_n})$ is not bigger than that of $\C$, which is 
$4\pi$, there is some $\beta_{n,j}$ in $Q_n^{-1}(\Gamma_j)$
of spherical length less than $C_M\eps/3$. Since 
$$J(b[s_2,1-s_2])=a[s_2,1-s_2]\subset\gamma^k\subset\B(0;M),$$ 
and 
$Q_n^{-1}$ converges to $J$ uniformly on $b[s_2,1-s_2]$, so if $n$ is big enough, then $Q_n^{-1}(b[s_2,1-s_2])
\subset \B(0;1.5M)$. Every curve in $\Gamma_j$ intersects $b[s_2,1-s_2]$, so $\beta_{n,j}\in Q_n^{-1}(\Gamma_j)$ 
intersects $Q_n^{-1}(b[s_2,1-s_2])\subset \mbox{\bf B}(0;1.5M)$. If $\beta_{n,j}\not\subset \B(0;2M)$, then there 
is a subarc $\gamma$ of $\beta_{n,j}$ that is contained in $\B(0;2M)$ and connects $\partial\B(0;1.5M)$ with $\partial\B(0;2M)$. So $L^\#(\gamma)\ge C_ML(\gamma)\ge C_MM/2$. This is impossible since 
$L^\#(\gamma)\le L^\#(\beta_{n,j})\le C_M\eps/3<C_MM/2$. 
Thus $\beta_{n,j}\subset \B(0;2M)$, and so $L(\beta_{n,j})\le L^\#(\beta_{n,j})/C_M<\eps/3$. 
Since $\beta_{n,j}$ has finite length, it is a crosscut in $D_n$ on $B_1^n$. 
Let $s_{n,0}$ be the biggest $s$ such that 
$Q_n^{-1}(b(s))\in\beta_{n,0}$, and $s_{n,1}$ the biggest $s$ such that $Q_n^{-1}(b(1-s))\in\beta_{n,1}$.
Then $s_{n,0},s_{n,1}\in[s_2,s_1]$. Let $\beta_{n,0}'$ and $\beta_{n,1}'$ denote any one component of 
$\beta_{n,0}\sem\{Q_n^{-1}(b(s_{n,0}))\}$ and $\beta_{n,1}\sem\{Q_n^{-1}(b(1-s_{n,1}))\}$, respectively.
Let $$\gamma_n^k:=Q_n^{-1}(b[s_{n,0},1-s_{n,1}])\cup\beta_{n,0}'\cup\beta_{n,1}'.$$ Then $\gamma_n^k$ is a
crosscut in $D_n$ on $B_1^n$. As $r_1<\eps$, the symmetric difference between $Q_n(\gamma_n^k)$ and 
$Q_0(\gamma^k)$ is contained in $\mbox{\bf B}(b(0);\eps)
\cup \mbox{\bf B}(b(1);\eps)$. Since $b[s_{n,0},1-s_{n,1}]$ is contained in $b[s_2,1-s_2]$, which
is a compact subset of $D_0$, so if $n$ is big enough, then
the Hausdorff distance between $Q_n^{-1}(b[s_{n,0},1-s_{n,1}])$ and $a[s_{n,0},1-s_{n,1}]$ is less than
$\eps/3$. Now the Hausdorff distance between $Q_n^{-1}(b[s_{n,0},1-s_{n,1}])$ and $\gamma_n^k$ is not
bigger than the bigger diameter of $\beta_{n,0}'$ and $\beta_{n,1}'$, which is less than $\eps/3$.
And the Hausdorff distance between $a[s_{n,0},1-s_{n,1}]$ and $\gamma^k$ is not bigger than the bigger diameter
of $a[0,s_{n,0}]$ and $a[1-s_{n,1},1]$, which is also less than $\eps/3$.
So the Hausdorff distance between $\gamma_n^k$ and $\gamma^k$ is less than $\eps$.
Now we proved that we can choose crosscuts $\gamma^{k}_n$ 
in $D_n$ on $B_1^n$ such that $\gamma^{k}_n$ converges to $\gamma^{k}$, and the symmetric difference of 
$Q_n(\gamma_n^{k})$ and $Q_0(\gamma^{k})$ converges to the two end points of $Q_0(\gamma^{k})$, 
respectively, both in the Hausdorff distance, as $n$ tends to infinity. 

\subsection{Constructing hooks that hold the boundary}

Now fix $k\ge 2$. We still parameterize $\lin{\gamma^k}$ and $\lin{Q_0(\gamma^k)}$ as the image of the function 
$a:[0,1]\to D\cup B_1^0$ and $b:[0,1]\to \A_{p_0}\cup \CC_0$, respectively, such that $b(t)=Q_0
(a(t))$, for $t\in (0,1)$. Let $\Omega^k$ denote the domain
bounded by $Q_0(\gamma^{k-1})$ and $Q_0(\gamma^{k+1})$ in $\A_{p_0}$. Then $\partial\Omega^k$ is composed of 
$Q_0(\gamma^{k-1})$, $Q_0(\gamma^{k+1})$, and two arcs on $\CC_0$. Let $\rho^k_0$ and $\rho^k_1$ denote
these two arcs such that $b(j)\in\rho^k_j$, $j=0,1$. If $n$ is big enough, 
from the convergence of $Q_n(\gamma_n^{k\pm 1})$ to $Q_0(\gamma^{k\pm 1})$, we have $\lin{Q_n(\gamma^{k-1}_n)}\cap
\lin{Q_n(\gamma^{k+1}_n)}=\emptyset$, and $U(Q_n(\gamma^{k+1}_n))\subset U(Q_n(\gamma^{k-1}_n))$.
Let $\Omega_n^k$ denote the domain bounded by $Q_n(\gamma_n^{k-1})$ and $Q_n(\gamma_n^{k+1})$ in $\A_{p_n}$.
Then the boundary of $\Omega_n^k$ is composed of $Q_n(\gamma_n^{k-1})$, $Q_0(\gamma_n^{k+1})$, and two disjoint 
arcs on $\CC_0$. If $n$ is big enough, then each of these two arcs contains one of $b(0)$ and $b(1)$. 
Let $\rho^k_{n,0}$ and $\rho^k_{n,1}$ denote these two arcs so that $b(j)\in\rho^k_{n,j}$, $j=0,1$. Now suppose $c:(-1,+1)\to \Omega^k$ is a crosscut in $\Omega^k$ with $c(\pm 1)
\in Q_0(\gamma^{k\pm 1})$. Then $c(-1,+1)$ divides $\Omega^k$ into two parts: $\Omega^k_0$ and $\Omega^k_1$,
so that $\rho^k_j\subset\partial\Omega^k_j$, $j=0,1$. If $n$ is big enough, then $c(\pm 1)\in Q_n(\gamma^{k\pm 1}_n)$,
and $c(-1,+1)\subset \Omega^k_n$. Thus $c(-1,+1)$ also divides $\Omega^k_n$ into two parts: $\Omega_{n,0}^k$
and $\Omega_{n,1}^k$, so that $\rho_{n,j}^k\subset\partial\Omega_{n,j}^k$. Let $\lambda_j$ 
($\lambda_{n,j}$, resp.) be the extremal distance
between $Q_0(\gamma^{k-1})$ ($Q_n(\gamma^{k-1}_n)$, resp.) and $Q_0(\gamma^{k+1})$ ($Q_n(\gamma^{k+1}_n)$, 
resp.) in $\Omega^k_j$ ($\Omega^k_{n,j}$, resp.), $j=0,1$. 
It is clear that $\lambda_{n,j}\to \lambda_j$ as $n\to\infty$, and $\lambda_j<\infty$.
Thus $\{\lambda_{n,j}\}$ is bounded by some $I_k>0$.

Since $\lin{\gamma^k}\cap\lin{\gamma^{k\pm 1}}=\emptyset$ and $\gamma^{k\pm 1}_n$ converges to $\gamma^{k\pm 1}$ 
in the Hausdorff distance, there is $d_k>0$ such that the distance between $\gamma^k$ and $\gamma_n^{k\pm 1}$
is greater than $d_k$, if $n$ is big enough.  For $x\in D_0$ and $r>0$, let $\til{\mbox{\bf B}}_0(x;r)$ 
and $\til{\mbox{\bf B}}_n(x;r)$ denote 
the connected component of $\mbox{\bf B}(x;r)\cap D_0$ and $\mbox{\bf B}(x;r)\cap D_n$, respectively,
 that contains $x$. 
Since $D_n\to D_0$, it is easy to prove that $\til{\mbox{\bf B}}_n(x;r)\to
\til{\mbox{\bf B}}_0(x;r)$. Let $e_k=d_k\exp(-2\pi I_k)$. Suppose $s_0\in(0,1)$ is such that the diameter 
of $a(0,s_0)$ is less than $e_k$. By the construction of $\gamma^k_n$, we have $\Omega^k_n\to\Omega^k$, 
so $Q_n^{-1}(\Omega^k_n)\to Q_0^{-1}(\Omega^k)$. Now $a(s_0)\in
\gamma^k\subset Q_0^{-1}(\Omega^k)$. Hence $a(s_0)\in Q_n^{-1}(\Omega^k_n)$ if $n$ is big enough.
Since the distance from $a(s_0)$ to $\gamma^{k\pm 1}_n$ is bigger than $d_k>e_k$, 
$\til{\mbox{\bf B}}
_n(a(s_0);e_k)$ is contained in $Q_n^{-1}(\Omega^k_n)$. We claim that $\til{\mbox{\bf B}}_n(a(s_0);e_k)
\subset Q_n^{-1}(\Omega^k_{n,0})$, if $n$ is big enough.

Since $a(0)\in\partial Q_0^{-1}(\Omega^k)$, $|a(0)-a(s_0)|<e_k$, and $Q_n^{-1}(\Omega^k_n)\to 
Q_0^{-1}(\Omega^k)$, so the distance from $a(s_0)$ to $\partial Q_n^{-1}(\Omega^k_n))$ is less than $e_k$, if 
$n$ is big enough. Now choose $z_n\in\partial Q_n^{-1}(\Omega^k_n)$ that is the nearest to $a(s_0)$. 
Then the line segment $[a(s_0),z_n)\subset\til{\mbox{\bf B}}_n(a(s_0);e_k)$. 
Hence $Q_n[a(s_0),z_n)$ is a simple curve in $\Omega^k_n$ such that $Q_n(z)$ tends to some $z_n'\in 
\partial \Omega^k_n$, as $z\in[a(s_0),z_n)$ and $z\to z_n$. Since $z_n\not\in\gamma^{k\pm 1}_n$, 
$z_n'\not\in Q_n(\gamma^{k\pm 1}_n)$. Thus $z_n'$ is on $\rho^k_{n,j}$ for some $j\in\{0,1\}$.
Since $Q_n(\til{\mbox{\bf B}}_n(a(s_0);e_k))\to Q_0(\til{\mbox{\bf B}}_0(a(s_0);e_k))
\ni b(s_0)$, and $b(s_0)\in\Omega^k_{n,0}$,
so if $n$ is big enough, $Q_n(\til{\mbox{\bf B}}_n(a(s_0);e_k))$ intersects $\Omega^k_{n,0}$.
For such $n$, if $z_n'\in\rho^k_{n,1}$, then 
all curves in $Q_n^{-1}(\Omega^k_{n,0})$ that go from
$\gamma^{k-1}_n$ to $\gamma^{k-1}_n$ will pass $\til{\mbox{\bf B}}_n(a(s_0);e_k)$. And so they all cross
some annulus centered at $a(s_0)$ with inner radius $e_k$ and outer radius greater than $d_k$. So the extremal
distance between $\gamma^{k-1}_n$ and $\gamma^{k+1}_n$ in $Q_n^{-1}(\Omega^k_{n,j})$ is greater than
$(\ln d_k-\ln e_k)/(2\pi)=I_k$. However, by conformal invariance, this extremal distance is equal to 
$\lambda_{n,j}$, which is not bigger than $I_k$ if $n$ is big enough. Thus 
$z_n'\in\rho^k_{n,0}$ for $n$ big enough. Similarly, $z_n'\in\rho^k_{n,0}$ and
$Q_n(\til{\mbox{\bf B}}_n(a(s_0);e_k))\cap\lin{\Omega^k_{n,1}\ne\emptyset}$ can not happen
at the same time when $n$ is big enough. So if $n$ is big enough, $Q_n(\til{\mbox{\bf B}}_n(a(s_0);e_k))$ is contained in $\Omega^k_{n,0}$. Similarly, we let 
$s_1\in(s_0,1)$ be such that the diameter of $a(s_1,1)$ is less than $e_k$, then 
$Q_n(\til{\B}_n(a(s_1);e_k))\subset\Omega^k_{n,1}$, if $n$ is big enough.

For $j=0,1$, $a(s_j)$ and $a(j)$ determine a square of side length $l_j=|a(j)-a(s_j)|$ with
vertices $v_{0,j}:=a(s_j)$, $v_{2,j}$, $v_{1,j}$, and $v_{3,j}$, in the clockwise order, so that 
$a(j)$ is on one middle line $[(v_{0,j}+v_{3,j})/2,(v_{1,j}+v_{2,j})/2]$. This square is contained in
$\lin{\mbox{\bf B}(a(s_j);\sqrt{2}l_j)}\subset \mbox{\bf B}(a(s_j);0.8 e_k)$, since $l_j<e_k/2$. 
And the union of line segments
$[v_{0,j},v_{1,j}]$, $[v_{1,j},v_{2,j}]$ and $[v_{2,j},v_{3,j}]$ surrounds $\mbox{\bf B}(a(j);l_j/8)$. 

For $j=0,1$, let $N_j$ be the $l_j/20$-neighborhood of $[v_{0,j},v_{1,j}]\cup[v_{1,j},v_{2,j}]\cup[v_{2,j},
v_{3,j}]$. Then $N_j\subset\B(a(s_j);e_k)$. Choose $q_j\in(0,l_j/30)$ such that 
$\lin{\B(a(s_j);q_j)}
\subset Q_0^{-1}(\Omega^k)$. For $m=0,1,2,3$, let $W_{m,j}=\lin{\B(v_{m,j};q_j)}$. When $n$ is big enough,
$W_{0,j}\subset Q_n^{-1}(\Omega_n^k)$, and $\B(a(j);l_j/30)$ intersects $\partial Q_n^{-1}(\Omega_n^k)$.
Suppose $\beta_j$ is a curve in $N_j$ which starts from $W_{0,j}$, and reaches $W_{1,j}$, $W_{2,j}$ and 
$W_{3,j}$ in the order. Then $\beta_j$ disconnects a subset of $\partial Q_n^{-1}(\Omega_n^k)$ from $\infty$,
if $n$ is big enough. Since $Q_n^{-1}(\Omega_n^k)$ is a simply connected domain,
$\beta_j$ hits $\partial Q_n^{-1}(\Omega_n^k)$. Let $\beta_j^n$ be the part of $\beta_j$ before hitting
$\partial Q_n^{-1}(\Omega_n^k)$. Then $\beta_j^n\subset\til{\B}_n(a(s_j);e_k)\subset Q_n^{-1}
(\Omega_{n,j}^k)$, if $n$ is big enough. So $Q_n(\beta_j^n)$ is a curve in $\Omega^k_{n,j}$ that tends to 
some point of $\partial\Omega^k_{n,j}$ at one end. This point is not on $Q_n(\gamma^{k\pm 1}_n)$, 
because the distance between $\gamma^k$ and $\gamma^{k+1}_n$ is greater than $e_k$.
Hence $\lin{Q_n(\beta_j^n)}$ intersects $\rho^k_{n,j}$.


Suppose $I$ is a closed ball in $Q_0^{-1}(\Omega^k)$. For $j=0,1$, let $\Pi_j$ be a subdomain of 
$Q_0^{-1}(\Omega^k)$ that contains $I\cup W_{0,j}$ such that $\lin{\Pi_j}$ is a compact subset of 
$Q_0^{-1}(\Omega^k)$. Then $\Pi_j$ is contained in $Q_n^{-1}(\Omega_n^k)$ for $n$ big enough. 
For $x\in\delta_n\Z^2\cap I$, let ${\cal A}_{n,j}^x$ be the set of lattice paths of $\delta_n\Z^2$ that start
from $x$, and hit $W_{0,j}$, $W_{1,j}$, $W_{2,j}$ and $W_{3,j}$ in the order before exiting $\Pi_j\cup N_j$.
We may view $\beta\in{\cal A}_{n,j}^x$ as a continuous curve. Let 
$\beta^{D_n}$ denote the part of $\beta\in{\cal A}_{n,j}^x$ before exiting $Q_n^{-1}(\Omega_n^k)$. Then 
$\beta^{D_n}$ can be viewed as a lattice path on $D^{\delta_n}$. We proved in the last paragraph that if $n$ is big 
enough, $\lin{Q_n(\beta^{D_n})}$ intersects $\rho^k_{n,j}$, for any $\beta\in{\cal A}_{n,j}^x$, 
$x\in\delta_n\Z^2\cap I$, $j=0,1$. Thus for any $\beta_0\in{\cal A}_{n,0}^x$ and 
$\beta_1\in{\cal A}_{n,1}^x$, $\beta_0^{D_n}\cup\beta_1^{D_n}$ 
disconnects $\gamma^{k-1}_n$ from $\gamma^{k+1}_n$ in $Q_n^{-1}(\Omega_n^k)$. 

\subsection{The behaviors of $g_0\circ J$ outside any neighborhood of $1$}

Let $P^x_{n,j}$ be the probability that a simple random walk on $\delta_n\Z^2$ started from $x$ belongs to 
${\cal A}_{n,j}^x$. By Lemma 3.8, if $n$ is big enough, then $P^x_{n,j}$ is greater than 
some $a_k>0$ for all $x\in\delta_n\Z^2\cap I$, $j=0,1$. We may also choose $n$ big enough 
such that
$V(D^{\delta_n})\cap I$ is non-empty, and $g_n(x)$ is less than some $b_k\in(0,\infty)$ for all 
$x\in\delta_n\Z^2\cap I$. We claim that if $n$ is big enough, then $g_n(x)\le\max\{b_k/a_k,1\}$
for every $x\in\delta_n\Z^2\cap(D_n\sem U(\gamma_n^{k-1}))$. Suppose for infinitely many $n$, there are $x_n\in\delta_n
\Z^2\cap D_n\sem U(\gamma_n^{k-1})$ such that $g_n(x_n)\ge M>\max\{b_k/a_k,1\}$. Since $g_n$ is discrete harmonic
on $\delta_n\Z^2\cap D_n$, and $g_n\le 1$ on the boundary vertices of $D_n$ except at $P(w_n)$,
 the tip point of 
$w_n$, so there is a lattice path $\beta_n$ in $D_n$ that goes from $x_n$ to $P(w_n)$ 
such that the value
of $g_n$ at each vertex of $\beta_n$ is not less than $M$. By the construction of 
$\gamma^{k+1}_n$, if $n$ is big enough, then $U(Q_n(\gamma^{k+1}_n))$ is some neighborhood of $1$ in $\A_{p_n}$, 
and so $U(\gamma^{k+1}_n)$ is some neighborhood of $P(w_n)$ in $D_n$. 
Thus $\beta_n$ intersects both $\gamma^{k-1}_n$ and $\gamma^{k+1}_n$.
Choose $v_0\in\delta_n\Z^2\cap I$. For every $\rho_{n,0}\in{\cal A}_{n,0}^{v_0}$ and
$\rho_{n,1}\in{\cal A}_{n,1}^{v_0}$, the path $\rho_{n,0}^{D_n}\cup\rho_{n,1}^{D_n}$ disconnects 
$\gamma^{k-1}_n$ from $\gamma^{k+1}_n$. Therefore $\rho_{n,0}^{D_n}\cup\rho_{n,1}^{D_n}$ intersects 
$\beta_n$. This implies that for some $j_n\in\{0,1\}$, for every $\rho\in{\cal A}_{n,j}^{v_0}$, we have
$\rho^{D_n}$ intersects $\beta_n$.
Thus the probability that a simple random walk on $\delta_{n}\Z^2$ started 
from $v_0$ hits $\beta_{n}$ before $\partial D_{n}$ is greater than $a_k$.
Let $\tau_n$ be the first time this random walk hits $\beta_n\cup\partial D_n$. Since $g_n$ is non-negative, bounded, 
and discrete harmonic on $\delta_n\Z^2\cap D_n$, so 
$g_n(v_0)=\mbox{\bf E}[g_n(\RW^x_{v_0}({\tau_n}))]\ge a_kM>b_k$, 
which is a contradiction. So the claim is proved.

By passing to a subsequence depending on $k$, we can now assume the following. 
$U(\gamma^{k+1}_n)$ is some neighborhood of $P(w_n)$ in $D_n$;
the value of $g_n$ on $\delta_n\Z^2\cap D_n
\sem U(\gamma_n^{k+1})$ is bounded by some $M_k\ge 1$;
$U(\gamma^{k+1}_n)\subset U(\gamma^{k}_n)\subset U(\gamma^{k-1}_n)$;  the spherical distance between 
$\gamma^{k}_n$ and $\gamma^{k-1}_n$ is greater than some $R_k>0$; and the (Euclidean) distance between
$\gamma^{k}_n$ and $\gamma^{k+1}_n$ is greater than $\delta_n$. Since the end points of $\gamma^{k}_n$ and $\gamma^{k-1}_n$
are on $B_1^n$, the spherical diameter of $B_1^n$ is at least $R_k$. Let $R$ be the spherical 
distance between $B_2$ and $\alpha_2$. Then the spherical
distance between $B_2$ and $B_1^n$ is at least $R$, as $\alpha_2$ disconnects $B_2$ from $B_1^n$.
Suppose $v\in V(D^{\delta_n})\cap D_n\sem U(\gamma_n^{k-1})$, and $dist^\#(v,B_1^n)=d<R/2$. 
Then $dist^\#(v,B_2)>R/2$. Let $\RW_v^n$ be a simple random walk on $\delta_n\Z^2$ started from $v$, 
and $\tau_n^{k}$ be the first time that $\RW_v^n$ leaves $D_n\sem U(\gamma_n^{k})$.
Then $\RW_v^n(\tau_n^{k})$ is either on $B_2$, or on $B_1^n$, or in $U(\gamma_n^{k})$. 
In the first case, $g_n(\RW_v^n(\tau_n^{k}))=1$, and $v$ should first exit $\B^\#(v,R/2)$
before hitting $B_2$. In the second and third cases, since 
$\RW_v^n(\tau_n^{k}-1)\in D_n\sem U(\gamma^{k}_n)$, and the Euclidean distance 
between $\gamma_n^k$
and $\gamma_n^{k+1}$ is greater than $\delta$ by construction, 
so $[\RW_v^n(\tau_n^{k}-1),\RW_v^n(\tau_n^{k})]$
does not intersect $\gamma^{k+1}_n$. Thus in the second case, $\RW_v^n(\tau_n^{k})\not=P(w_n)$,
and so $g_n(\RW_v^n(\tau_n^{k}))=0$. In the third case, $\RW_v^n(\tau_n^{k})\in D_n\sem
U(\gamma^{k+1}_n)$, so $g_n(\RW_v^n(\tau_n^{k}))\le M_k$; and $\RW_v^n$
first uses some edge that intersects $\gamma^{k-1}_n$, then uses some edge that intersects
$\gamma^{k}_n$ at time $\tau_n^{k}$. So the spherical diameter of $\RW_v^n[0,\tau_n^{k}]$ is
at least $R_k$. This implies that $\RW_v^n$ should first exit $\B^\#(v;R_k/2)$ before
hitting $U(\gamma^{k}_n)$. Let $R_k'=\min\{R/2,R_k/2\}$, then by Lemma 3.7, 
$$\mbox{\bf P}[\RW_v^n(\tau_n^{k})\not\in B_1^n]\le C_0((\delta_n+d)/R_k')^{C_1},$$
for some absolute constants $C_0,C_1>0$. 
So we have $g_n(v)\le M_kC_0((\delta_n+d)/R_k')^{C_1}$.

Suppose $z\in D_0\sem U(\gamma^{k-1})\sem\gamma^{k-1}$, and $dist^\#(z,B_1^0)=d<R/4$. Choose 
$r\in (0,d/2)$ such that $\B^\#(z,r)$ is bounded and 
$\lin{\B^\#(z;r)}\subset D_0\sem U(\gamma^{k-1})\sem\gamma^{k-1}$.
If $n$ is big enough, then $\lin{\B^\#(z;r)}\subset D_n\sem U(\gamma_n^{k-1})$, and
the spherical distance from every $v\in\B^\#(z;r)$ to $B_1^n$ is less than $2d<R/2$.
Thus $$g_n(v)\le M_kC_0((\delta_n+2d)/R_k')^{C_1}\mbox{, }\mbox{ }\forall v\in\delta_n\Z^2\cap\B^\#(z;r).$$
Since $g_0$ is the limit of $g_n$, $g_0(z)\le M_kC_0(2d/R)^{C_1}$. Thus 
for every $k\ge 2$, $g_0(z)\to 0$,
as $z\in D_0\sem U(\gamma^{k-1})\sem\gamma^{k-1}$, and $z\to B_1$ in the spherical metric, and 
so $g_0\circ J(z)\to 0$ as $z\in\A_{p_0}\sem U(Q_0(\gamma^{k-1}))$, and $z\to
\CC_0$. 
Since $U(Q_0(\gamma^k))$, $k\in\N$, forms a neighborhood basis of $1$ in $\A_{p_0}$, 
so for any $r>0$, $g_0\circ J(z)\to 0$ if $z\in \A_{p_0}\sem\B(1,r)$ and $z\to\CC_0$. 
This is what we need at the end of 5.1. $\Box$
\vskip 4mm

{\bf Acknowledgment.} This work was proceeded under the instruction of Professor Nikolai 
Makarov, who let the author be interested in this subject, and gave many valuable comments
on this paper.

\vskip 4mm

-----------------------------------------

Department of Mathematics

Mail code: 253-37

California Institute of Technology

Pasadena, CA 91125, USA

dapeng@its.caltech.edu

\end{document}